\documentclass[10pt]{amsart}
\usepackage{amssymb,amsmath,amsthm,amscd,MnSymbol}
\usepackage{epsfig, comment,enumerate, color}

\usepackage[all]{xy}

\newtheorem{thm}{Theorem}
\newtheorem{cor}[thm]{Corollary}
\newtheorem{lem}[thm]{Lemma}
\newtheorem{prop}[thm]{Proposition}
\theoremstyle{definition}
\newtheorem{defn}[thm]{Definition}

\newtheorem{setup}[thm]{Setup}

\newtheorem{rk}[thm]{Remark}

\newtheorem{theorem-question}[thm]{Theorem-Question}

\newenvironment{exafont}{\begin{bf}}{\end{bf}}

\newcommand{\ttB}{\mathtt{B}}
\newcommand{\bfc}{\mathbf{c}}
\newcommand{\bfa}{\mathbf{a}}\newcommand{\bfb}{\mathbf{b}}
 \newcommand{\bft}{\mathbf{t}}
 \newcommand{\ut}{\underline{\mathbf{t}}}

\newcommand{\um}{\underline{\mathbf{m}}}

\newcommand{\bbO}{\mathbb{O}}\newcommand{\fO}{\mathfrak{O}}

\newcommand{\cA}{\mathcal{A}}
\newcommand{\bbH}{\mathbb{H}}\newcommand{\bbHH}{\mathbb{HH}}
\newcommand{\ffHH}{\mathfrak{HH}}
\newcommand{\cB}{\mathcal{B}}

\newcommand{\ZZ}{\mathbb{Z}}\newcommand{\hh}{\mathbf{hh}}

\newcommand{\bbT}{\mathbb{T}}

\newcommand{\x}{\mathbf{x}}\newcommand{\bfm}{\mathbf{m}}
\newcommand{\y}{\mathbf{y}}
\newcommand{\GL}{\mathrm{GL}}
\newcommand{\HTct}{\boldsymbol{\Lambda}}\newcommand{\HHTct}{\boldsymbol{\Pi}}

\DeclareMathOperator{\HH}{HH}
\DeclareMathOperator{\End}{End}
\DeclareMathOperator{\Hom}{Hom}
\DeclareMathOperator{\RHom}{RHom}

\DeclareMathOperator{\Ext}{Ext}

\DeclareMathOperator{\ml}{-mod}

\DeclareMathOperator{\mlr}{-mod-}

\DeclareMathOperator{\rbigr}{bigr_{jk}-}
\DeclareMathOperator{\lbigr}{-bigr_{jk}}
\DeclareMathOperator{\bibigr}{-bigr_{jk}-}

\DeclareMathOperator{\perf}{-perf}\DeclareMathOperator{\rperf}{perf-}
\DeclareMathOperator{\op}{op}

\title{Hochschild cohomology of polynomial representations of $\GL_2$}

\author{Vanessa Miemietz and Will Turner}

\begin{document}

\begin{abstract}
We compute the Hochschild cohomology algebras of Ringel-self-dual blocks of polynomial representations of $\GL_2$ over an algebraically closed field of
characteristic $p>2$, that is, of any block whose number of simple modules is a power of $p$. These algebras are finite-dimensional and we provide an explicit description of their bases and multiplications.
\end{abstract}

\maketitle \vspace{-1cm}
\tableofcontents
\setlength{\parskip}{8pt}

\section{Introduction.}

Hochschild cohomology is a basic invariant which associates to a finite dimensional algebra $A$  a super-commutative algebra 
$\HH(A) = \Ext^\bullet_{A \mlr A}(A)$.
The algebra $\HH(A)$ can be thought of as the derived centre of the algebra $A$, 
given as it is by the formula $\HH(A) = H^\bullet\End_{A \mlr A}(\tilde{A})$, 
where $\tilde{A}$ is a projective resolution of $A$ in the category $A \mlr A$ of $A$-$A$-bimodules; 
to see the analogy compare with the formula $Z(A) =\End_{A \mlr A}(A)$ for the classical centre $Z(A)$ of a unital algebra $A$.
If $M$ is any $A$-module, then the natural algebra homomorphism $Z(A) \rightarrow \Hom_A(M,M)$ extends to a 
natural algebra homomorphism $\HH(A) \rightarrow \Ext^\bullet_A(M,M)$.

Like other algebras obtained by taking derived endomorphisms,
Hochschild cohomology and its variants can be endowed with additional structures, which have been the source of diverse interest:
the most basic such is known as the Gerstenhaber bracket \cite{Ger}.
But even without further decoration, the algebra $\HH(A)$ has proved difficult to compute in 
specific examples, and its behaviour difficult to predict. 
One delicacy is the issue of finite generation of $\HH(A)$ which is not guaranteed for a finite dimensional algebra $A$, 
even modulo the ideal of nilpotent elements \cite{Xu, Snashall}; 
yet there are finite dimensional self-injective algebras whose Hochschild cohomology 
is not merely finitely generated but finite dimensional \cite{BGMS}.

The subject of this article is the computation of $\HH$ in a basic example arising in the representation theory of algebraic groups.
We examine the Hochschild cohomology of polynomial representations of the algebraic group $G = \GL_2(F)$, 
where $F$ is an algebraically closed field of characteristic $p$. 
Indeed, we compute the Hochschild cohomology of any Ringel self-dual block of polynomial representations of $G$ for $p>2$, which by \cite[Theorem 27]{EH} are precisely those blocks with $p^l$ simple modules for $l \in \mathbb{N}$.
The algebras describing these blocks increase in complexity as $l$ increases, 
but we are nevertheless able to develop sufficiently sharp homological tools to achieve the calculation of their $\HH$ algebras. Their Hochschild cohomology algebras, for which we give explicit bases and multiplications, turn out not only to be finitely generated, but indeed finite-dimensional.
We apply a theory of algebraic operators ($2$-functors) on certain $2$-categories which underlies the representation theory of $G$ \cite{MT2}, \cite{MT3}.
We also use the theory of quasi-hereditary algebras \cite{CPS}, the theory of Koszul duality \cite{BGS}, the formalism of differential graded algebras and their derived categories \cite{Ke}, a theory of homological duality for algebraic operators,  explicit analysis of certain bimodules associated with a well-known quasi-hereditary algebra $\bfc$ and its homological duals, and a formalism of algebras with a polytopal basis.

\section{The answer.}

All algebras considered in this article will be $F$-algebras.
Suppose $\Gamma = \bigoplus_{i,j,k \in \mathbb{Z}} \Gamma^{ijk}$ is a $\mathbb{Z}$-trigraded algebra.
We have a combinatorial operator $\fO_\Gamma$ which acts on the collection of $\mathbb{Z}$-bigraded algebras $\Sigma$ after the formula
$$\fO_\Gamma(\Sigma)^{ik} = \bigoplus_{j, k_1+k_2=k} \Gamma^{ijk_1} \otimes_F \Sigma^{jk_2},$$
where we take the super tensor product with respect to the $k$-grading.

Let $p>2$.
In the main body of the paper we define an $ijk$-graded algebra $\HHTct$ with an explicit, 
canonically defined basis $\mathcal{B}_{\HHTct}$. 
A complete description of the algebra $\HHTct$, its basis, and its product, is given in Section \ref{spadesuitsection}.

There is a natural algebra homomorphism $F \leftarrow \HHTct$ which is a splitting of the map sending $1$ to the identity in $\HHTct$.
This lifts to a morphism of operators $\fO_F \leftarrow \fO_{\HHTct}$, which means that we obtain an algebra homomorphism $\fO_F\Sigma \leftarrow \fO_{\HHTct}\Sigma$ for every $\Sigma$. Since $\fO_F^2 = \fO_F$ we obtain a sequence of operators
$$\fO_F \leftarrow \fO_F \fO_{\HHTct} \leftarrow \fO_F \fO_{\HHTct}^2 \leftarrow \fO_F \fO_{\HHTct}^3 \leftarrow ...$$
We define $\hh_l$ to be the Hochschild cohomology of a block of polynomial representations of $G$ with $p^l$ simple modules
and establish the following:

\begin{thm} \label{main}
For any $l>0$, the algebra $\hh_l$ is isomorphic to $\fO_F \fO_{\HHTct}^l (F[z,z^{-1}])$.
\end{thm}

\begin{rk} For every $l$ the algebra $\hh_l$ inherits an explicit basis from $\HHTct$ with an explicit product as described in Corollary \ref{totalbasis}.
\end{rk}

{\bf Acknowledgement.} The first author acknowledges support from ERC grant PERG07-GA-2010-268109. We would also like to thank the referee for an extremely thorough and helpful report.

\section{Guidebook.}

The proof of Theorem \ref{main} passes through a number of counties of diverse character; 
here we briefly describe some of these. 
The algebras we are interested in are not Koszul algebras; nevertheless, they are closely related to certain Koszul algebras and
we make use of some pretty generalities concerning the Hochschild cohomology of Koszul algebras; 
in Section \ref{hhkoszulsection} we give an account of these.
In Section \ref{oldsection} we introduce certain algebraic 
operators and gather together some facts about these that we have established in previous papers.
In Section \ref{algopHH} we describe an interaction of these operators with Hochschild cohomology and Koszul duality.
In Section \ref{GL2section} we recall from another paper \cite{MT2} how special examples of our algebraic operators can be used to 
describe the  polynomial representation theory of $\GL_2(F)$.
In Section \ref{computation} we show that this description of the polynomial representation theory of $\GL_2(F)$
via algebraic operators along with the Section \ref{algopHH} 
analysis of the behaviour of Hochschild cohomology under such algebraic operators can be used to describe the 
Hochschild cohomology for the algebras relevant to $\GL_2(F)$ in terms of an algebraic operator $\fO_{\ffHH(\HTct)}$; 
here $\HTct = \bbH\bbT_{\Omega}(\ut^!))$ is the homology tensor algebra over a  certain Koszul algebra $\Omega$ 
of a certain pair of dg $\Omega$-$\Omega$ bimodules $\ut^!$,
and $\ffHH$ is the operator that sends a graded algebra $X = \bigoplus_i X^i$ to a graded algebra $\bigoplus_i \HH(X^0, X^i)$.
In Section \ref{clubsuitsection} we give a combinatorial description of the algebra $\HTct$ via certain bimodules;
to do this we invoke a study of the negative part $\HTct^-$ of $\HTct$ made in a previous article \cite{MT3}, 
and Serre duality for $\Omega$.
In Section \ref{explicitbimod} we perform a detailed combinatorial analysis of the Hochschild cohomology of 
certain bimodules appearing in the algebra $\HTct$.
A fact emerging here is that a certain quotient $\Theta$ of $\Omega$,
commonly known as the preprojective algebra of type $A$, possesses an involution $\sigma$ such that 
$$\Theta^\sigma \cong \Theta^*, \qquad \qquad \HH(\Omega, \Theta^\sigma) \cong \HH(\Omega, \Theta)^*;$$
the first of these formulas asserts the well known self-injectivity of $\Theta$, 
but the second asserts something similar holds under $\HH(\Omega,-)$.
In Section \ref{spadesuitsection} we use the analysis of the preceding section to give a combinatorial description of 
$\HHTct = \bbHH(\HTct)$ in terms of certain bimodules and maps between them.
Finally in Section \ref{final} we reach our destination, and give a proof of Theorem \ref{main} as well as a monomial basis for the algebras we construct.

\section{Hochschild cohomology of Koszul algebras.} \label{hhkoszulsection}
\subsection{Grading conventions.}\label{grad}
In order to fix our notations, we will now give a brief introduction to dg algebras and modules, which will be the main objects of study in this paper.
A differential graded vector space is a $\mathbb{Z}$-graded vector space $V = \oplus_k V^k$  with a graded endomorphism $d$ of degree $1$.
We write $|v|$ for the degree of a homogeneous element of $V$. We will always assume all $V^k$ to be finite-dimensional.
We assume $d$ can act both on the left and the right of $V$, with the convention $d(v) = (-1)^{|v|}(v)d$.
A differential graded algebra is a $\mathbb{Z}$-graded algebra $A = \oplus_k A^k$ with a differential $d$ such that
$$d(ab) = d(a).b + (-1)^{|a|}a.d(b),$$
or equivalently
$$(ab)d = a.(b)d + (-1)^{|b|}(a)d.b.$$
If $A$ is a differential graded algebra then a differential graded left $A$-module is a graded left $A$-module $M$ with differential $d$ such that
$$d(a.m) = d(a).m + (-1)^{|a|}a.d(m);$$
a differential graded right $A$-module is a graded right $A$-module $M$ with differential $d$ such that
$$d(m.a) = d(m).a + (-1)^{|m|}m.d(a).$$
If $A$ and $B$ are dg algebras then a dg $A$-$B$-bimodule is a graded $A$-$B$-bimodule with a differential
which is both a left dg $A$-module and a right dg $B$-module.

If $_AM_B$ and $_BN_C$ are dg bimodules where $A$, $B$, and $C$ are dg algebras, then $M \otimes_B N$ is a dg $A$-$C$-bimodule
with differential $$d(m \otimes n) = d(m) \otimes n + (-1)^{|m|} m \otimes d(n).$$

Speaking about morphisms of dg algebras and dg (bi-)modules we mean homogeneous morphisms.
However, if $_AM_B$ and $_AN_C$ are dg bimodules where $A$, $B$, and $C$ are dg algebras, then $\Hom_A(M, N)$, the $k$-graded vector space whose $k$-degree $m$-part consists of all $A$-module morphisms $f:M\to N$ such that $f(M^\bullet) \subseteq M^{\bullet+m}$. This  is a dg $B$-$C$-bimodule with differential $$d(\phi) = d\circ\phi - (-1)^{|\phi|} \phi\circ d.$$

If $_AM$ is a left dg $A$-module,
then $\End_A(M)$ is a differential graded algebra which acts on the right of $M$, giving $M$ the structure of an $A$-$\End_A(M)$-bimodule,
the differential on $\End_A(M)$ being given by $(\phi)d = \phi\circ d - (-1)^{|\phi|} d \circ\phi$.
If $M_B$ is a right dg $B$-module,
then $\End_B(M)$ is a differential graded algebra which acts on the left of $M$, giving $M$ the structure of an $\End_B(M)$-$B$-bimodule,
the differential on $\End_B(M)$ being given by $d(\phi) =  d\circ\phi - (-1)^{|\phi|} \phi \circ d$.

A differential bi- (tri-)graded vector space is a vector space $V$  with a $\mathbb{Z}^2$- respectively $\mathbb{Z}^3$-grading whose coordinates we denote by $(j,k)$ respectively $(i,j,k)$ and an endomorphism $d$ of degree $(0,0,1)$, i.e. $d$ is homogeneous with respect to the $i,j$-gradings and has degree $1$ in the $k$-grading, which we will also call the homological grading.
We denote by $\langle \cdot \rangle$ a shift by $1$ in the $j$-grading, meaning $(V
\langle n \rangle)^{j} = V^{j-n}$. Since we will often identify dg modules and complexes, we will stick to the complex convention of $[\cdot]$ being a shift to the left, i.e.\ $V[n]^k = V^{k+n}$. Altogether $$(V\langle n \rangle [m])^{ijk} = V^{i,j-n,k+m}.$$
All definitions above can be extended to the differential bi- (tri-)graded setting, defining differential bi- (tri-)graded algebras, differential bi- (tri-)graded (left and right) $A$-modules as well as bi- (tri-)graded  $A$-$B$-bimodules as bi- (tri-)graded algebras resp.~modules resp.~bimodules which are differential graded algebras resp.~modules resp.~bimodules with respect to the $k$-grading, i.e. with respect to an endomorphism of degree $(0,0,1)$. Speaking about morphisms of differential bi- (tri-)graded  algebras and differential bi- (tri-)graded  (bi-)modules we mean homogeneous morphisms with respect to all gradings.
Similarly to the above, homomorphism spaces taken between $A$-modules (rather than differential (bi-) trigraded $A$-modules) will carry a differential bi- (tri-)grading.

For a dg algebra $A$, we denote by $D_{dg}(A)$ the dg derived category of $A$, whose objects are (left) dg $A$-modules and where morphisms are given by the localisation of the class of dg module morphisms with respect to the class which are quasi-isomorphisms (see \cite[Section 3.1, 3.2]{Ke}). We denote by $A\perf$ and $\rperf A$ the  categories of left resp.~right perfect dg $A$-modules.

We let 
$\mathbb{H}$ denote the cohomology functor, which takes a differential
$k$-graded complex $C$ to the $k$-graded vector space $\mathbb{H} C
= H^{\bullet} C$.

\subsection{Hochschild cohomology  of Koszul algebras.} Koszul duality was introduced by Beilinson, Ginzburg and Soergel \cite{BGS} and generalised to dg algebras by Keller \cite[Section 10]{Ke1}. The conventions we follow are given in \cite[Appendix B]{MT4}, and also summarised below. 

\begin{setup}\label{koszulconv}
Throughout this section, $A$ denotes a finite-dimensional Koszul algebra. It is hence in particular a quadratic $j$-graded algebra of the form 
$A = \bbT_{A^0}(A^1)/R$, with relations $R \subset A^1 \otimes_{A^0} A^1$, and we write
$A^! = \bbT_{A^{0}}((A^!)^{-1})/R^!$ for its quadratic dual (which is then also Koszul), where the $A^0$-$A^0$ bimodules $A^1$ and $(A^!)^{-1}$,
and the short exact sequences of $A^0$-$A^0$-bimodules
$$0 \rightarrow R \rightarrow A^1 \otimes_{A^0} A^1 \rightarrow A^2 \rightarrow 0$$
$$0 \leftarrow (A^!)^{-2} \leftarrow (A^!)^{-1}\otimes_{A^0} (A^!)^{-1} \leftarrow R^{!} \leftarrow 0,$$
are duals of each other.
We insist $A$ is generated in $j$-degrees $0$ and $1$, and $A^!$ is generated in $j$-degrees $0$ and $-1$. We assume that $A^0$ is isomorphic to a direct product of a number of copies of $F$, and denote by $e_s$ the idempotent corresponding to the $s$th copy.
\end{setup}

Following \cite[Proposition 2.2.4.1]{LH} and \cite[Section 4.7]{Ke2}, the Koszul resolution is given by $A\otimes_\tau (A^!)^*\otimes_\tau A$ where $\tau$ is the canonical twisting cochain given by the composition $$(A^!)^* \to A^1 \to A$$ of the inclusion by the projection.
The complex $A\otimes_\tau (A^!)^*\otimes_\tau A$ is isomorphic as a complex to $B:= A\otimes_{A^0} (A^!)^*\otimes_{A^0} A$, with differential
$$\alpha\otimes \varphi\otimes \alpha' \mapsto  \sum_{\rho \in \ttB^1} ((-1)^{|\alpha|}\alpha \rho \otimes \rho^*\varphi \otimes \alpha' - (-1)^{|\varphi|+|\alpha|+|\rho^*|} \alpha\otimes \varphi \rho^* \otimes \rho \alpha'),$$
where $\ttB^1$ is any basis of the free $A^0$-$A^0$-bimodule $A^1$ (cf. \cite[page 1119]{Ne}). 

It follows from \cite[Theorem 6.3]{Ne} that there is an isomorphism of dg algebras
$$ \Hom_{A^0\otimes A^{0\op}} (A^{!*}, A) \to \Hom_{A\otimes A^{\op}} (A\otimes_{A^0} A^{!*}\otimes_{A^0} A, A\otimes_{A^0} A^{!*}\otimes_{A^0} A) $$
given by
$$f\mapsto \left( \alpha\otimes \varphi\otimes \alpha'\mapsto (-1)^{|f|(|\alpha|+|\varphi_1|)} \alpha\otimes \varphi_{(1)}\otimes f(\varphi_{(2)})\alpha'\right)$$
where the algebra structure on $\Hom_{A^0\otimes A^{0\op}} (A^{!*}, A) $ is induced by the comultiplication $\Delta :A^{!*} \to A^{!*}\otimes_{A^0} A^{!*}$ on $A^{!*}$ and we write $\Delta(\varphi)=\varphi_{(1)} \otimes \varphi_{(2)}$. Note that the original source considers tensor products and hom spaces over $F$, but the results readily generalise to our setup.

Let now $X, Y$ be differential $jk$-graded $A$-$A$-bimodules.
It then follows similarly that  the morphism
\begin{equation}\label{negroniso}
 \Hom_{A^0\otimes A^{0\op}} (A^{!*}, X) \to \Hom_{A\otimes A^{\op}} (A\otimes_{A^0} A^{!*}\otimes_{A^0} A, A\otimes_{A^0} A^{!*}\otimes_{A^0} X) 
 \end{equation}
given by
$$f\mapsto \left( \alpha\otimes \varphi\otimes \alpha'\mapsto (-1)^{|f|(|\alpha|+|\varphi_1|)} \alpha\otimes \varphi_{(1)}\otimes f(\varphi_{(2)})\alpha'\right)$$
translates
the product 
$$\Hom_{A^0\otimes A^{0\op}} (A^{!*}, X)\otimes \Hom_{A^0\otimes A^{0\op}} (A^{!*}, Y) \to \Hom_{A^0\otimes A^{0\op}} (A^{!*}, X\otimes_A Y) $$
induced by comultiplication on $A^{!*}$ into the cup product $$\HH(A,X) \otimes \HH(A,Y)\to  \HH(A,X\otimes_A Y) $$ after taking homology.

\begin{lem} \label{adjunctionsequence} In the situation of Setup \ref{koszulconv}, and for $X$ a differential $jk$-bigraded $A$-$A$-bimodule, we have isomorphisms of $jk$-graded vector spaces,
\begin{equation*}\begin{split}
\Hom_{A \otimes A^{\op}}(B,X) &\cong \Hom_{A^0 \otimes A^{0 \op}}(A^0, A^! \otimes_{A^0} X) \\
&\cong \bigoplus_{s,t} e_s A^! e_t \otimes_F e_t X e_s \\
&\cong \bigoplus_{s,t} (e_s \otimes e_s)(A^! \otimes_F X^{ \op}) (e_t \otimes e_t).
\end{split}\end{equation*}
Explicitly, the isomorphism
$$\bigoplus_{s,t} e_s A^! e_t \otimes_F e_t X e_s \to \Hom_{A \otimes A^{\op}}(B,X)$$
is given by 
$$a\otimes x \mapsto \left(\chi_{a\otimes x}\colon\;\alpha\otimes\varphi\otimes \alpha'\mapsto (-1)^{|\varphi||x|+(|a|+|x|)|\alpha|}\alpha\varphi(a)x\alpha' \right).$$
\end{lem}
\proof The second and third isomorphisms hold by definition. The first holds by a sequence of adjunctions:
\begin{equation*}\begin{split}
\Hom_{A \otimes A^{\op}}(A \otimes_{A^0} A^{!*} \otimes_{A^0} A, X)
&\cong \Hom_{A^0 \otimes A^{0 \op}}(A^{!*}, X)\\
&\cong \Hom_{A^0 \otimes A^{0 \op}}(A^0, \Hom_{A^0}(A^{!*}, X))\\
&\cong \Hom_{A^0 \otimes A^{0 \op}}(A^0, A^! \otimes_{A^0} X).\end{split}\end{equation*}

Tracing these adjunctions, we obtain the desired isomorphism $$\bigoplus_{s,t} e_s A^! e_t \otimes_F e_t X e_s \to \Hom_{A \otimes A^{\op}}(B,X)$$ as the composition of the isomorphism \begin{equation}\label{firstiso}\bigoplus_{s,t} e_s A^! e_t \otimes_F e_t X e_s \to\Hom_{A^0 \otimes A^{0 \op}}(A^{!*}, X)\end{equation}
given by
$$ a\otimes x\mapsto \left(\varphi \mapsto (-1)^{|\varphi||x|} \varphi(a)x\right)$$
and the isomorphism 
$$\Hom_{A^0 \otimes A^{0 \op}}(A^{!*}, X)\to  \Hom_{A \otimes A^{\op}}(A \otimes_{A^0} A^{!*} \otimes_{A^0} A, X)$$
given by
$$f \mapsto \left( \alpha\otimes \varphi\otimes \alpha'\mapsto (-1)^{|f||\alpha|} \alpha f(\varphi)\alpha'\right),$$
which is the composition of the morphism in \eqref{negroniso} and the natural projection $A \otimes_{A^0} A^{!*} \otimes_{A^0} X\to X$. 
\endproof

\begin{setup}\label{gradconv}
In addition to keeping the conventions from Setup  \ref{koszulconv}, we now further assume we are in one of the following cases:

\begin{enumerate}
\item\label{conv1}$A$ is $jk$-graded such that $A^{1\bullet} = A^{10}$, and $(A^!)^{-1\bullet}=(A^!)^{-11}$. This implies in particular that $A$ is concentrated in $k$-degree $0$, and $A^!$ is concentrated in non-negative $k$-degrees.
\item\label{conv2} $A$ is $jk$-graded such that $A^{1\bullet} = A^{11}$, and $(A^!)^{-1\bullet}=(A^!)^{-10}$. This implies in particular that $A^!$ is concentrated in $k$-degree $0$, and $A$ is concentrated in non-negative $k$-degrees.
\end{enumerate}
\end{setup}

In particular, assuming Setup \ref{gradconv} \eqref{conv1}, the differential on $B=A \otimes_{A^0} A^{!*} \otimes_{A^0} A$ specialises to  
$$\alpha\otimes \varphi\otimes \alpha' \mapsto  \sum_{\rho \in \ttB^1} (\alpha \rho \otimes \rho^*\varphi \otimes \alpha' + (-1)^{|\varphi|} \alpha\otimes \varphi \rho^* \otimes \rho \alpha')$$
and, assuming Setup \ref{gradconv} \eqref{conv2}, the differential on $B$ specialises to specialises to
$$\alpha\otimes \varphi\otimes \alpha'\mapsto (-1)^{|\alpha|}\sum_{\rho \in \ttB^1} (\alpha \rho \otimes \rho^*\alpha \otimes \alpha' -  \alpha\otimes \varphi \rho^* \otimes \rho \alpha').$$

\begin{thm} \label{bimodulehochschild} 
Assume we are in the situation of Setup \ref{gradconv}\eqref{conv2} and let $X$ be a differential $jk$-bigraded $A$-$A$-bimodule.
Then we  have an isomorphism
$$\HH(A,X) \cong \bbH \left(\bigoplus_{s,t}  e_sA^!e_t \otimes e_t \bbH X e_s\right)$$
where the differential on $\bigoplus_{s,t} e_sA^!e_t \otimes e_t \bbH X e_s$ is given by
$$a\otimes x \mapsto - \sum_{\rho \in \ttB^1} (a\rho^*\otimes \rho x- (-1)^{|x|} \rho^*a\otimes x\rho).$$
\end{thm}

\proof 
In Setup \ref{gradconv}\eqref{conv2}, the $jk$-graded vector space isomorphism 
$$\bigoplus_{s,t} e_s A^! e_t \otimes_F e_t X e_s \to  \Hom_{A \otimes A^{\op}}(A \otimes_{A^0} A^{!*} \otimes_{A^0} A, X)$$
 from  Lemma \ref{adjunctionsequence} is now given by
\begin{equation}\label{chidef}
a\otimes x \mapsto 
\left(\chi_{a\otimes x}\colon \;  \alpha\otimes\varphi\otimes \alpha'\mapsto (-1)^{|\alpha||x|}\alpha\varphi(a)x\alpha' \right),\end{equation}
which, wanting this to be an isomorphism of $jk$-graded differential bimodules, forces the differential on $\bigoplus_{s,t} e_s A^! e_t \otimes_F e_t X e_s$ to be given by
$$ d \colon a\otimes x \mapsto (1\otimes d_X)(a\otimes x) - \sum_{\rho \in \ttB^1} (a\rho^*\otimes \rho x -(-1)^{|x|} \rho^*a\otimes x\rho )$$
where $d_X$ again denotes the differential on $X$. 

Indeed, we compute the differential of $\chi_{a\otimes x}$ which is the map 
\begin{equation*}
\begin{split}
\alpha\otimes\varphi\otimes \alpha'&\mapsto (-1)^{|\alpha||x|}(-1)^{|\alpha|}\alpha\varphi(a)d_X(x)\alpha' \\
&-(-1)^{|x|}\chi_{a\otimes x}\left( (-1)^{|\alpha|}\sum_{\rho \in \ttB^1} (\alpha \rho \otimes \rho^*\varphi \otimes \alpha' -  \alpha\otimes \varphi \rho^* \otimes \rho \alpha') \right)\\
&= (-1)^{|\alpha|(|x|+1)}\alpha\varphi(a)d_X(x)\alpha'\\
&-(-1)^{|\alpha|+|x|} (-1)^{(|\alpha|+1)|x|}\sum_{\rho \in \ttB^1}
\alpha\rho(\rho^*\varphi)(a)x\alpha' \\
&+(-1)^{|\alpha|+|x|}(-1)^{|\alpha||x|}\sum_{\rho \in \ttB^1}\alpha(\varphi\rho^*)(a)x\rho\alpha' \\
&= (-1)^{|\alpha|(|x|+1)}\alpha\varphi(a)d_X(x)\alpha'\\
&-(-1)^{|\alpha|(|x|+1)}\left( \sum_{\rho \in \ttB^1}
\alpha\rho\varphi(a\rho^*)x\alpha \right) \\
&-(-1)^{(|\alpha|+1)(|x|+1)}\left( \sum_{\rho \in \ttB^1}\alpha\varphi(\rho^*a)x\rho\alpha \right).
\end{split}
\end{equation*}

On the other hand, with the prescribed differential $d$, 
$$\chi_{d(a\otimes x)}= \chi_{a\otimes d_X(x)} - \sum_{\rho \in \ttB^1} \left(\chi_{a\rho^*\otimes \rho x}-(-1)^{|x|} \chi_{\rho^*a\otimes x\rho }\right),$$
which is the map
\begin{equation*}
\begin{split}
\alpha\otimes\varphi\otimes \alpha'&\mapsto (-1)^{|\alpha|(|x|+1)}\alpha \varphi(a)d_x(x)\alpha'\\
&-(-1)^{|\alpha|(|x|+1)}\sum_{\rho \in \ttB^1}\alpha\varphi(a\rho^*)\rho x \alpha'\\
&+(-1)^{|x|}(-1)^{|\alpha|(|x|+1)} \sum_{\rho \in \ttB^1} \alpha\varphi(\rho^*a) x\rho\alpha'\\
&= (-1)^{|\alpha|(|x|+1)}\alpha \varphi(a)d_X(x)\alpha'\\
&-(-1)^{|\alpha|(|x|+1)}\sum_{\rho \in \ttB^1}\alpha\varphi(a\rho^*)\rho x \alpha'\\
&-(-1)^{(|\alpha|+1)(|x|+1)} \sum_{\rho \in \ttB^1} \alpha\varphi(\rho^*a) x\rho\alpha'\\
\end{split}
\end{equation*}
which equals the expression for the differential of $\chi_{a\otimes x}(\alpha\otimes \varphi\otimes \alpha')$ term by term.

We write $d= 1\otimes d_X +\tilde d$ for the differential on $\bigoplus_{s,t} e_s A^! e_t \otimes_F e_t X e_s$. We are interested in the homology of this complex.
Notice that the map (which is in itself a differential) $a \otimes x \mapsto   a\otimes d_X(x)$ has $j$-degree $0$ on each tensor factor, and the remaining part of the differential has $j$-degree $-1$ on the first and $j$-degree $1$ on the second tensor factor. 
We now consider the spectral sequence induced by the radical filtration of $A^!$. 
Then it follows immediately from the definition that $d_0 = 1\otimes d_X$, $d_1 = \tilde d$ and $d_l= 0$ for all $l\geq 2$ and consequently $$\bbH(\bigoplus_{s,t} e_sA^!e_t \otimes e_t X e_s)\cong \bbH(\bigoplus_{s,t} e_sA^!e_t \otimes e_t \bbH X e_s)$$ where the differential on the latter complex is precisely given by $\tilde d$.  
\endproof

We now consider the cup product in Hochschild cohomology.

\begin{prop}\label{cup} Assume we are in the situatiom of Setup \ref{gradconv}\eqref{conv2} and 
let $X$ and $Y$ be differential $jk$-bigraded $A$-$A$-bimodules. Under the isomorphism in Theorem \ref{bimodulehochschild}, the cup product 
$$\HH(A,X) \otimes \HH(A,Y)\to  \HH(A,X\otimes Y)$$
is translated to
$$\begin{array}{rcl}
\left(\bigoplus_{s,t}  e_sA^!e_t \otimes e_t \bbH X e_s\right) \otimes\left(\bigoplus_{s,t}  e_sA^!e_t \otimes e_t \bbH Y e_s\right)&\to & \left(\bigoplus_{s,t}  e_sA^!e_t \otimes e_t \bbH( X\otimes_A Y) e_s\right)\\
(a\otimes x) \otimes (b\otimes y)&\mapsto& ba\otimes xy.
\end{array}$$
\end{prop}

\proof
Recall from \eqref{negroniso} the isomorphism of dg vector spaces
$$\mathbf{T}\colon \Hom_{A^0\otimes A^{0\op}} (A^{!*}, X) \to \Hom_{A\otimes A^{\op}} (A\otimes_{A^0} A^{!*}\otimes_{A^0} A, A\otimes_{A^0} A^{!*}\otimes_{A^0} X) $$
given by
$$f\mapsto \left( \alpha\otimes \varphi\otimes \alpha'\mapsto (-1)^{|f|(|\alpha|+|\varphi_1|)} \alpha\otimes \varphi_{(1)}\otimes f(\varphi_{(2)})\alpha'\right).$$
Notice that composing this with the natural projection $A\otimes_{A^0} A^{!*}\otimes_{A^0} X \to X$, the only term in the sum that survives in $\mathbf{T}(f)(\alpha\otimes \varphi\otimes \alpha')$ is $(-1)^{|f|(|\alpha|)} \alpha\otimes 1\otimes f(\varphi)\alpha'$ which maps to $(-1)^{|f|(|\alpha|)} \alpha f(\varphi)\alpha'$. 
Hence, if $f$ is the image of $a\otimes x$ under the isomorphism \eqref{firstiso}, we precisely obtain our $\chi_{a\otimes x}$ from \eqref{chidef} above. 

The product 
$$\Hom_{A^0\otimes A^{0\op}} (A^{!*}, X)\otimes \Hom_{A^0\otimes A^{0\op}} (A^{!*}, Y) \to \Hom_{A^0\otimes A^{0\op}} (A^{!*}, X\otimes_A Y) $$
induced by comultiplication on $A^{!*}$ translates into the cup product $\HH(A,X) \to \HH(A,Y)\to  \HH(A,X\otimes Y) $ after taking homology.
We thus consider the translation of the product 
$$\Hom_{A^0\otimes A^{0\op}} (A^{!*}, X)\otimes \Hom_{A^0\otimes A^{0\op}} (A^{!*}, Y) \to \Hom_{A^0\otimes A^{0\op}} (A^{!*}, X\otimes_A Y) $$ induced by comultiplication on $A^{!*}$
into a product 
$$\left(A^!\otimes_{A^0\otimes A^{0\op}}X\right) \otimes\left(A^!\otimes_{A^0\otimes A^{0\op}}Y\right)\to A^!\otimes_{A^0\otimes A^{0\op}}\left(X\otimes_AY\right).$$
Denoting the image of $a\otimes x$ under the isomorphism \eqref{firstiso} by $\xi_{a\otimes x}\colon (\varphi\mapsto \varphi(a)x)$
(using $|\varphi|=0$),
we see that from
$(\xi_{a\otimes x}\cdot\xi_{b\otimes y})(\varphi)= \varphi_{(1)}(a)\varphi_{(2)}(b)x\otimes y$
it follows that the product of $a\otimes x$ and $b\otimes y$, being the preimage of $(\xi_{a\otimes x}\cdot\xi_{b\otimes y})$ is $ba\otimes (x\otimes y)$, from the formula $\Delta(\varphi)(a\otimes b) =\varphi(ba)$. 

Since the splitting of the differential on $A^!\otimes_{A^0\otimes A^{0\op}}X$ as $d=d_X+\tilde d$ is compatible with the tensor product, the isomorphism obtained in Theorem \ref{bimodulehochschild} is compatible with the cup product in homology and we obtain the desired multiplication formula.\endproof

\section{Some old things.} \label{oldsection}

Here we gather an assortment of notions and facts we have established in previous articles.
More details can be found in those articles \cite{MT2}, \cite{MT3}.

\subsection{The $2$-category $\mathcal{T}$}
Let $\mathcal{T}$ denote the
collection of  pairs $(A, M)$ where $A$ is a differential $k$-graded algebra and $M$
is a  differential $k$-graded $A$-$A$-bimodule.

The collection $\mathcal{T}$ in fact forms the set of objects of a $2$-category:  $1$-morphisms  from $(A,M)$ to $(B,N)$ are given by a pair $(S, \phi_S)$, consisting of a differential (bi-)graded $A$-$B$-bimodule ${}_AS_B$ and a quasi-isomorphism
$$\phi_S: S \otimes_BN\to M\otimes_AS; $$
$2$-morphisms from $(S, \phi_S)$ to $(T, \phi_T)$ are given by homomorphisms of differential (bi-)graded $A$-$B$-bimodules $f: S \to T$ such that the diagram
$$\xymatrix{
S \otimes_BN\ar^{\phi_S}[r]\ar^{ f\otimes id }[d]&M \otimes_A S \ar^{id \otimes f}[d] \\
T \otimes_BN\ar^{\phi_T}[r]                      &M \otimes_A T                     
}$$
commutes.

\begin{defn}
We define a {\bf Rickard object} of  $\mathcal{T}$ to be an object
$(A,M)$ of $\mathcal{T}$, where ${}_AM_A$ is perfect as a left dg $A$-module and as a right dg $A$-module, the natural morphism of dg algebras $A\to \RHom_A(M,M)$ is a quasi-isomorphism, there is a quasi-isomorphism $A \to \bbH A$, and $\bbH A$ is a finite-dimensional algebra of finite global dimension. We call a Rickard object $(A,M)$ a {\bf classical Rickard object} if $A$ has zero differential, and ${}_AM_A$ is projective on both sides.
\end{defn}

\begin{defn}
We define a {\bf $j$-graded
object} of $\mathcal{T}$ to be an object $(a,m)$ of $\mathcal{T}$,
where $a= \bigoplus a^{jk}$ is a differential bigraded algebra, and $m
= \bigoplus m^{jk}$ a differential bigraded $a$-$a$-bimodule, and $a^{j\bullet} = m^{j\bullet} = 0$ for $j<0$. 
\end{defn}

\subsection{The operator $\mathbb{O}$.} Let $(\bfa,\bfm) $ be a $j$-graded object of $\mathcal{T}$. We define
$$\mathbb{O}_{\bfa,\bfm} \circlearrowright \mathcal{T}$$
to be the operator given by
$$\mathbb{O}_{\bfa,\bfm}(A,M) = (\bfa(A,M), \bfm(A,M)),$$
where
$$\alpha(A,M) =(\alpha^0 \otimes A)\oplus (\bigoplus_{j > 0} \alpha^j \otimes M^{\otimes_A j})$$ for $\alpha \in \{\bfa,\bfm\}.$ 
The algebra structure on $\bigoplus \bfa^{jk} \otimes_F M^{\otimes_A j}$ is the restriction of the algebra structure on the tensor
product of algebras $\bfa \otimes \bbT_{A}(M)$, where $ \bbT_{A}(M)$ denotes the tensor algebra of $M$ over $A$.
The $k$-grading and differential on the complex $\bigoplus \bfa^{jk} \otimes M^{\otimes_A j}$ are defined to be the total $k$-grading and total differential on the tensor product of complexes. The bimodule structure, grading and differential on $\bigoplus \bfm^{jk} \otimes M^{\otimes_A j}$ are defined likewise. 

We remark that this extends to a $2$-endofunctor of $\mathcal{T}$ (cf. \cite[Lemma 9]{MT2}).
\begin{lem} \label{tensor}\cite[Lemma 14]{MT3}
Let $\bfa, \bfb, \bfc$ be a differential bigraded algebras, ${}_\bfb \x_\bfa$ and ${}_\bfa \y_\bfc$
differential bigraded modules, all concentrated in nonnegative $j$-degrees. Let  $(A, M)$ be an object of $\mathcal{T}$.
Then
$$\x(A,M) \otimes_{\bfa(A,M)} \y(A,M) \cong (\x \otimes_\bfa \y)(A,M)$$ as differential bigraded $\bfb(A,M)$-$\bfc(A,M)$-bimodules.
\end{lem}

Given a differential bigraded $\bfa$-module $\x$, with components in
positive and negative $j$-degrees, we define $\x(A,M)$ to be the
$\bfa(A,M)$-module given by
$$\x(A,M) = (\bigoplus_{j<0} x^{j\bullet} \otimes (M^{-1})^{\otimes_A -j}) \oplus (x^{0\bullet} \otimes A) \oplus
(\bigoplus_{j>0} x^{j\bullet} \otimes M^{\otimes_A j}),$$
where $M^{-1}:=\Hom_A(M,A)$.

\begin{lem} \label{homo}(cf.\cite[Lemma 15]{MT3})
Let $\bfc$ be a differential bigraded algebra with $\bfc^0$ a product of copies of $F$. Let $\x$ and $\y$ are differential bigraded
$\bfc$-modules, all concentrated in nonnegative $j$-degrees, and let $(A,M)$ be a Rickard
object of $\mathcal{T}$. Then we have a quasi-isomorphism of differential bigraded $(\bfc^0\otimes A)\otimes (\bfc^0\otimes A)^{op}$-modules
$$\Hom_\bfc(\x,\y)(A,M) \rightarrow \Hom_{\bfc(A,M)}(\x(A,M), \y(A,M)).$$
\end{lem}

\proof
This was proved in  \cite[Proof of
Theorem 13]{MT2} (though it was stated only as a quasi-isomorphism of differential bigraded vector spaces there), but for the convenience of the reader we recall the construction.

For notational simplicity, we write $M^j:=M^{\otimes_A j}$ for $j\geq0$ and $M^j:=M^{\otimes_A -j}$ for $j\leq 0$, and also  write $t_1 \cdots t_j :=t_1\otimes\cdots \otimes t_j\in M^j$ for $j>0$.

We construct a map
$$\beta: \Hom_{\bfc}(\x,\y)(A,M) \rightarrow \Hom_{\bfc(A,M)}(\x(A,M), \y(A,M)).$$ 

For an element $f_j \otimes t_1 \cdots t_j$, or $f_j \otimes h$, 
where $f_j \in \Hom_\bfc^{(j)}(\x,\y)$, and $t_1 \cdots t_j \in M^j$ if $j \geq 0$, 
and where $h \in \Hom_A(M^{-j},A)$ if $j<0$, we define

$$\beta(f_j \otimes t_1 \cdots t_j)= 
(x_k \otimes t_1'\cdots t_k' \mapsto f_j(x_k) \otimes  t_1'\cdots t_k't_1 \cdots t_j),$$
$$\beta(f_j \otimes h)=\left(x_k \otimes t_1'\cdots t_k'\mapsto
\left\{ \begin{array}{ll}
 0& \hbox{if } -j > k\\
 f_j(x_k) \otimes t_1'\cdots t_{k-(-j)}'h(t_{k-(-j)+1}'\dots t_k')  & \hbox{if } -j \leq k
\end{array}\right)\right.,$$ 
where we work with the convention that for all $A$-$A$-bimodules $M$ and $N$ such that $N$ is finitely generated and projective, we identify $\Hom_A(N,A)\otimes_A \Hom_A(M,A)$ with $\Hom_A(M\otimes_A N,A)$ via the map sending any $g\otimes f$ to the morphism from $M\otimes_A N$ to $A$ given by $m\otimes n \mapsto g(n)f(m)$.
Using an explicit quasi-isomorphism
\begin{equation*}
\begin{split}\Hom_A(M^{i-1}, A\otimes_A M^{j-1}) &\rightarrow \Hom_A(M^{i-1}, \Hom_A(M,M) \otimes_A M^{j-1}) \\
&  \rightarrow  \Hom_A(M^{i},M^{j})\end{split}\end{equation*} 
given by
\begin{equation*}
\begin{split}( t_1 \cdots t_{i-1} \mapsto 1 \otimes f(t_1 \cdots t_{i-1}))&\mapsto ( t_1 \cdots t_{i-1} \mapsto (t_i \mapsto t_i) \otimes f(t_1 \cdots t_{i-1}))\\
& \mapsto (t'_1 \cdots t'_{i}  \mapsto  t'_1 f(t_2' \cdots t'_i)).\end{split}\end{equation*} 
iteratively, and applying the technical result \cite[Lemma 16]{MT2}, which computes the space $\Hom_{\bfc(A,M)}(\x(A,M), \y(A,M))$ as a $(\bfc^0\otimes A)$-$(\bfc^0\otimes A)$-bimodule,  one sees that $\beta$ is a quasi-isomorphism  of $(\bfc^0\otimes A)$-$(\bfc^0\otimes A)$-bimodules.
\endproof

\subsection{ The operator $\fO$.} We now recall the definition of the operator $\fO$ from \cite{MT3}. Let $\Gamma = \bigoplus \Gamma^{ijk}$ be a
differential trigraded algebra. We have an operator
$$\fO_\Gamma \circlearrowright \{ \Sigma | \mbox{ $\Sigma = \bigoplus \Sigma^{jk}$ a differential bigraded algebra } \}$$
given by
\begin{equation}\label{fO}\fO_\Gamma(\Sigma)^{ik} = \bigoplus_{j, k_1+k_2=k} \Gamma^{ijk_1} \otimes \Sigma^{jk_2}.\end{equation}
The algebra structure and differential are obtained by restricting
the algebra structure and differential from $\Gamma \otimes \Sigma$.
If we forget the differential and the $k$-grading, the operator
$\fO_\Gamma$ is identical to the operator $\fO_\Gamma$ defined in the
introduction.

Suppose we are given $(\bfa_i,\bfm_i)$ for $1 \leq i \leq n$, and $(A,M)$.
Let us define $(A_i,M_i)$ recursively via $(A_i,M_i) =
\bbO_{\bfa_i,\bfm_i}(A_{i-1},M_{i-1})$ and $(A_0,M_0) = (A,M)$.

\begin{lem} \label{comparetwo}\cite[Lemma 20]{MT3}
\begin{enumerate}[(i)]
 \item \label{9i} We have an algebra isomorphism
$$\bbT_{A_n}(M_n) \cong \fO_{\bbT_{\bfa_n}(\bfm_n)}
...\fO_{\bbT_{\bfa_1}(\bfm_1)}(\bbT_A(M)).$$
\item\label{9ii} We have an isomorphism
of objects of $\mathcal{T}$
\begin{equation*}\begin{split}\bbO_{\bfa_n,\bfm_n} &\cdots\bbO_{\bfa_1,\bfm_1}(A,M) \\
&\cong (\fO_{\bbT_{\bfa_1}(\bfm_1)}...\fO_{\bbT_{\bfa_n}(\bfm_n)}(\bbT_A(M))^{0\diamond \bullet},
\fO_{\bbT_{\bfa_1}(\bfm_1)}...\fO_{\bbT_{\bfa_n}(\bfm_n)}(\bbT_A(M))^{1\diamond \bullet}).
\end{split}\end{equation*}
\end{enumerate}
\end{lem}

\section{Algebraic operators and Hochschild cohomology.}\label{algopHH}

Given a Rickard object $(\bfa,\bfm)$ in $\mathcal{T}$, we define $\bbHH(\bfa,\bfm):=\bigoplus_{i\in \ZZ}\HH(\bfa,\bfm^{\otimes_\bfa i})$, where for $i<0$, we understand $\bfm^{\otimes_\bfa i}$ as $(\bfm^{-1})^{\otimes_\bfa{-i}}$ for $\bfm^{-1} = \Hom_{\bfa}(\bfm, \bfa)$ the bimodule adjoint to $\bfm$.

\begin{lem}\label{algstruc}
Let $(\bfa,\bfm)$ be a $j$-graded classical Rickard object in $\mathcal{T}$.
Then the space $\bbHH(\bfa,\bfm)$ has the structure of an $ijk$-trigraded associative algebra.
\end{lem}

\proof
Let $\tilde \bfa$ be a projective $\bfa$-$\bfa$-bimodule resolution of $\bfa$ and note that this implies that $\tilde \bfm^i:= \tilde \bfa\otimes_{\bfa} \bfm^{\otimes_\bfa i}$ is a complex of projective bimodules quasi-isomorphic to $\bfm^{\otimes_\bfa i}$. 
In this case, we have natural isomorphisms
\begin{equation*}\begin{split}
\bbH \Hom_{\bfa \otimes \bfa^{\op}}(\tilde \bfm^{h},\tilde \bfm^{i+h})&\cong \bbH \Hom_{\bfa \otimes \bfa^{\op}}(\tilde \bfa, \Hom_{\bfa^{op}}(\bfm^{\otimes_\bfa h},\tilde \bfm^{i+h}))\\
&\cong \bbH \Hom_{\bfa \otimes \bfa^{\op}}(\tilde \bfa, \tilde \bfm^i\otimes _\bfa \Hom_{\bfa^{op}}(\bfm^{\otimes_\bfa h},\bfm^{\otimes_\bfa h}))\\
&\cong  \bbH \Hom_{\bfa \otimes \bfa^{\op}}(\tilde \bfa, \tilde \bfm^{i})\\
&=\HH(\bfa,\bfm^{\otimes_\bfa i})  .\end{split}\end{equation*}
for any $h\in\ZZ$, where the first isomorphism is just adjunction, the second relies on $\bfm^{\otimes_\bfa h}$ being finitely generated projective as a right $\bfa$-module, and the third comes from $(\bfa,\bfm)$ being Rickard.
Identifying $\HH(\bfa,\bfm^{\otimes_\bfa i})$  with $\bbH \Hom_{\bfa \otimes \bfa^{\op}}(\tilde \bfm^{h},\tilde \bfm^{h+i})$ via this isomorphism gives us an associative product
$$\HH(\bfa,\bfm^{\otimes_\bfa i}) \otimes \HH(\bfa,\bfm^{\otimes_\bfa h}) \rightarrow \HH(\bfa,\bfm^{\otimes_\bfa h+i})$$
that equips $\mathbb{HH}(\bfa, \bfm)$ with the structure of an $ijk$-graded algebra. Note that this is the algebra structure obtained from the general definition of the cup product.
\endproof

\begin{thm} \label{iteration}
Let $(\bfa,\bfm)$ be a $j$-graded object in $\mathcal{T}$ with $\bfa$ concentrated in $k$-degree zero and Koszul, and let $(A,M)$ be a Rickard object in   $\mathcal{T}$.
Then we have
$$\bbHH \bbO_{\bfa, \bfm}( A,M) \cong \fO_{\bbHH(\bfa, \bfm)}(\mathbb{HH}( A,M))$$
as ijk-graded algebras.
\end{thm}
\proof
Since $\bfa$ is Koszul, we have  a projective $\bfa$-$\bfa$-bimodule resolution of $\bfa$ given by $$(\bfa \otimes_{\bfa^0} \bfa^{!*} \otimes_{\bfa^0} \bfa) \twoheadrightarrow \bfa$$ by Section \ref{hhkoszulsection}.
Let now $\tilde A$ be a projective $A$-$A$-bimodule resolution of $A$, and in analogy to the above, $\tilde M^{i}= \tilde A \otimes_A M^{\otimes_A i}$ a corresponding $A$-$A$-bimodule resolution of  $M^{\otimes_A i}$ for any $i\in \ZZ$. 

We now set $ \bfa^{!*}(\widetilde{A, M}):=(\bfa^{!*})^0\otimes \tilde A \oplus \bigoplus_{j> 0 } (\bfa^{!*})^j\otimes \tilde M^j$. (Note that since $\bfa^!$ is negatively $j$-graded, its dual is again positively $j$-graded.)
We claim that $$\bfa(A,M) \otimes_{\bfa^0(A,M)} \bfa^{!*}(\widetilde{A, M}) \otimes_{\bfa^0(A,M)} \bfa(A,M) \twoheadrightarrow \bfa(A,M)$$ is a projective $\bfa(A,M)$-$\bfa(A,M)$ bimodule resolution.

Indeed, as $\bfa^{!*}$ is projective over $\bfa^0 \otimes \bfa^{0\op}$ 
and $\tilde M^j$ is projective over $A \otimes A^{\op}$ for every $j$, we have that $ \bfa^{!*}(\widetilde{A, M})$ is projective over $$\bfa^0\otimes A\otimes (\bfa^0\otimes A)^{\op} = \bfa^0(A,M)\otimes \bfa^0(A,M)^{\op}.$$ Furthermore, $\bfa(A,M)$ is projective over $\bfa^0(A,M)=\bfa^0 \otimes A$ on both sides so $\bfa(A,M) \otimes \bfa(A,M)^{\op}$ is projective in $\bfa^0(A,M) \otimes \bfa^0(A,M)^{\op} \ml$. Therefore, the induced  module $\bfa(A,M) \otimes_{\bfa^0(A,M)} \bfa^{!*}(\widetilde{A,M}) \otimes_{\bfa^0(A,M)} \bfa(A,M)$ is projective in $\bfa(A,M) \otimes \bfa(A,M)^{\op} \ml$. By construction, it is quasi-isomorphic to the bimodule
$\bfa(A,M) \otimes_{\bfa^0(A,M)} \bfa^{!*}(A, M) \otimes_{\bfa^0(A,M)} \bfa(A,M)$, which by Lemma \ref{tensor} is quasi-isomorphic to 
  $(\bfa \otimes_{\bfa^0} \bfa^{!*} \otimes_{\bfa^0} \bfa) (A,M)$ and hence to $\bfa(A,M)$.

Now, setting $$\widetilde{\bfa(A,M)}:=\bfa(A,M) \otimes_{\bfa^0(A,M)} \bfa^{!*}(\widetilde{A, M}) \otimes_{\bfa^0(A,M)} \bfa(A,M),$$ we have isomorphisms
\begin{equation*}\begin{split}
&\HH(\bfa(A,M), \bfm(A,M)^{\otimes_{\bfa(A,M)} i}) \\&\cong \bbH\Hom_{\bfa(A,M) \otimes \bfa(A,M)^{\op}}(\widetilde{\bfa(A,M)},\bfm(A,M)^{\otimes_{\bfa(A,M)} i})\\
& \cong \bbH\Hom_{\bfa^0(A,M) \otimes \bfa^0(A,M)^{\op}}( 
\bfa^{!*}(\widetilde{A,M}) ,\bfm(A,M)^{\otimes_{\bfa(A,M)} i}) \\
& = \bbH\Hom_{\bfa^0\otimes A\otimes (\bfa^0\otimes A)^{\op}}( \bfa^{!*}(\widetilde{A,M}) ,\bfm(A,M)^{\otimes_{\bfa(A,M)} i}) \end{split}
\end{equation*}
by projectivity of $\widetilde{\bfa(A,M)}$ as an $\bfa(A,M) \otimes \bfa(A,M)^{\op}$-module and adjunctions. 

Next, notice that, as an $\bfa^0\otimes A$-$\bfa^0\otimes A$-bimodule,  $\bfa^{!*}(\widetilde{A,M}) \cong (\bfa^0\otimes \tilde A) \otimes_{ \bfa^0\otimes A} \bfa^{!*}(A,M)$, so we can use adjunction again and obtain
\begin{equation}\label{eq0}\begin{split}
\bbH&\Hom_{\bfa^0\otimes A\otimes (\bfa^0\otimes A)^{\op}}( \bfa^{!*}(\widetilde{A,M}) ,\bfm(A,M)^{\otimes_{\bfa(A,M)} i})\\&\cong\bbH\Hom_{\bfa^0\otimes A\otimes (\bfa^0\otimes A)^{\op}}( (\bfa^0\otimes \tilde A) \otimes_{ \bfa^0\otimes A} \bfa^{!*}(A,M) ,\bfm(A,M)^{\otimes_{\bfa(A,M)} i})\\
&\cong \bbH\Hom_{\bfa^0\otimes A\otimes (\bfa^0\otimes A)^{\op}}(\bfa^0\otimes \tilde A, \Hom_{(\bfa^0\otimes A)^{op}}(\bfa^{!*}(A,M) ,\bfm(A,M)^{\otimes_{\bfa(A,M)} i}))
 \end{split}
\end{equation}

We now claim that $\bfm(A,M)^{\otimes_{\bfa(A,M)} i}$ is quasi-isomorphic to $\bfm^{\otimes_\bfa i}(A,M)$ as $\bfa^0\otimes A$-$\bfa^0\otimes A$-bimodules. This follows directly from Lemma \ref{tensor} for $i>0$. For $i<0$, we obtain
\begin{equation*}\begin{split}
\bfm(A,M)^{\otimes_{\bfa(A,M)} i} &= \Hom_{\bfa(A,M)}(\bfm(A,M), \bfa(A,M))^{\otimes_{\bfa(A,M)} -i}\\
&\cong  \Hom_{\bfa(A,M)}(\bfm(A,M)^{\otimes_{\bfa(A,M)}- i}, \bfa(A,M))\\
&\cong \Hom_{\bfa(A,M)}(\bfm^{\otimes_{\bfa} -i}(A,M), \bfa(A,M))\\
& \leftarrow^{qim}  \Hom_{\bfa}(\bfm^{\otimes_{\bfa} -i}, \bfa)(A,M)\\
&\cong \Hom_{\bfa}(\bfm, \bfa)^{\otimes_{\bfa} -i}(A,M)\\
&= \bfm^{\otimes_{\bfa} i}(A,M)
 \end{split}
\end{equation*}
where the first  isomorphism comes from iterated adjunction and the fact that for $(\bfa,\bfm)$ Rickard, $\bfm(A,M)$ is again perfect as a left dg $A$-module and as a right dg $A$-module, the second isomorphism is Lemma \ref{tensor}, the quasi-isomorphism is Lemma \ref{homo} and the final isomorphism again uses  iterated adjunction and the fact that $\bfm$ is perfect as a left dg $\bfa$-module and as a right dg $\bfa$-module.

Using this, we have a quasi-isomorphism of $\bfa^0\otimes A$-$\bfa^0\otimes A$-bimodules
\begin{equation*}\begin{split}
\Hom_{(\bfa^0\otimes A)^{\op}}&(\bfa^{!*}(A,M) ,\bfm(A,M)^{\otimes_{\bfa(A,M)} i})\\& \leftarrow^{qim} \Hom_{(\bfa^0\otimes A)^{\op}}(\bfa^{!*}(A,M) ,\bfm^{\otimes_\bfa i}(A,M))\\
&\cong  \Hom_{(\bfa^0(A,M))^{\op}}(\bfa^{!*}(A,M) ,\bfm^{\otimes_\bfa i}(A,M))\\
& \leftarrow^{qim} \Hom_{(\bfa^0)^{\op}}(\bfa^{!*} ,\bfm^{\otimes_\bfa i})(A,M)\\
&\cong (\bfm^{\otimes_\bfa i}\otimes_{\bfa^0}\bfa^{!})(A,M).
 \end{split}
\end{equation*}

Putting this back into \eqref{eq0}, we obtain 
\begin{equation*}\begin{split}
\bbH\Hom_{\bfa^0\otimes A\otimes (\bfa^0\otimes A)^{\op}}&(\bfa^0\otimes \tilde A, \Hom_{(\bfa^0\otimes A)^{\op}}(\bfa^{!*}(A,M) ,\bfm(A,M)^{\otimes_{\bfa(A,M)} i}))\\
&\cong \bbH\Hom_{\bfa^0\otimes A\otimes (\bfa^0\otimes A)^{\op}}(\bfa^0\otimes \tilde A, (\bfm^{\otimes_\bfa i}\otimes_{\bfa^0}\bfa^{!})(A,M))\\
&\cong  \bbH\Hom_{\bfa^0\otimes A\otimes (\bfa^0\otimes A)^{\op}}(\bfa^0\otimes \tilde A, \bigoplus_{j}(\bfm^{\otimes_\bfa i}\otimes_{\bfa^0}\bfa^{!})^j\otimes M^{\otimes_A j})\\
&\cong \bbH \bigoplus_{j}\Hom_{\bfa^0\otimes \bfa^{0\op}}(\bfa^0, (\bfm^{\otimes_\bfa i}\otimes_{\bfa^0}\bfa^{!})^j)\otimes \Hom_{A\otimes A^{\op}}(\tilde A, M^{\otimes_A j})\\
&\cong \bbH \bigoplus_{j}\Hom_{\bfa^0\otimes \bfa^{0\op}}(\bfa^0, (\bfm^{\otimes_\bfa i}\otimes_{\bfa^0}\bfa^{!}))^j\otimes \Hom_{A\otimes A^{\op}}(\tilde A, M^{\otimes_A j})\\
&\cong \bbH \bigoplus_{j}\Hom_{\bfa\otimes \bfa^{\op}}(\bfa \otimes_{\bfa^0} \bfa^{!*} \otimes_{\bfa^0} \bfa,\bfm^{\otimes_\bfa i})^j \otimes \Hom_{A\otimes A^{\op}}(\tilde A, M^{\otimes_A j})\\
&=\bbH\Hom_{\bfa\otimes \bfa^{\op}}(\bfa \otimes_{\bfa^0} \bfa^{!*} \otimes_{\bfa^0} \bfa,\bfm^{\otimes_\bfa i})(\bbHH(A,M))\\
&=  \bbO_{\bbHH(\bfa,\bfm)}(\bbHH(A,M))^{i\bullet\diamond}.
 \end{split}
\end{equation*}
Summing over all $i$, we obtain the desired isomorphism as $ijk$-graded vector spaces.
In order to check that this is an isomorphism of $ijk$-graded algebras, we set 
$$\tilde\bfa:=\bfa \otimes_{\bfa^0} \bfa^{!*} \otimes_{\bfa^0}\bfa, \qquad \tilde\bfm^h:=\tilde\bfa  \otimes_{\bfa} \bfm^{\otimes_\bfa h}$$ and 
$$\widetilde{\bfm(A,M)}^h:=\widetilde{A(A,M)} \otimes_{\bfa(A,M)} \bfm(A,M)^{\otimes_{\bfa(A,M)} h}.$$ One then constructs a similar isomorphism as above for 
\begin{equation*}\begin{split}\bbH\Hom_{\bfa(A,M)\otimes \bfa(A,M)^{\op}}&(\widetilde{\bfm(A,M)}^h, \widetilde{\bfm(A,M)}^{h+i})\\&\cong \bbH\Hom_{\bfa\otimes \bfa^{\op}}(\tilde{\bfm}^h, \tilde{\bfm}^{h+i}) (\bbHH(A,M)) \end{split}
\end{equation*}
and checks that due to the naturality of all constructions the isomorphisms in the diagram
$$\xymatrix{  \bbH\Hom_{\bfa(A,M)\otimes \bfa(A,M)^{\op}}(\widetilde{\bfm(A,M)}^h, \widetilde{\bfm(A,M)}^{h+i})\ar@{<->}[r]\ar@{<->}[d]&\ar@{<->}[d] \bbH\Hom_{\bfa\otimes \bfa^{\op}}(\tilde{\bfm}^h, \tilde{\bfm}^{h+i}) (\bbHH(A,M))\\
\bbH\Hom_{\bfa(A,M)\otimes \bfa(A,M)^{\op}}(\widetilde{\bfa(A,M)},\bfm(A,M)^{\otimes_{\bfa(A,M)} i})\ar@{<->}[r]\ & \bbH\Hom_{\bfa\otimes \bfa^{\op}}(\tilde\bfa , \bfm^{\otimes_\bfa i}) (\bbHH(A,M))
 }$$
commute. \endproof

Observe that for a $j$-graded Rickard object $(\bfa,\bfm)$ in $\mathcal{T}$, the differential trigraded $\bfa$-$\bfa$-bimodule $\bigoplus_{i\in \ZZ} \bfm^{\otimes_\bfa i}$ (where $\bfm^{\otimes_\bfa 0}:=\bfa$), obtains the structure of an associative $ijk$-trigraded algebra when passing to homology, coming from the natural (quasi-) isomorphisms $\bfm\otimes_\bfa\bfm^{-1} \overset{\sim}{\to} \bfa $ of evaluation and $\bfm^{ -1}\otimes_\bfa\bfm\overset{\sim}{\to} \End_\bfa(\bfm) \overset{qim}{\leftarrow}\bfa.$ We denote this $ijk$-trigraded algebra by $\bbH\bbT_\bfa(\um)$ where $\um$ stands for $(\bfm, \bfm^{-1})$.
We define the $ijk$-graded vector space
$$\ffHH(\bbH\bbT_\bfa(\um)):= \bigoplus_{i\in \ZZ} \HH(\bfa, \bbH(\bfm^{\otimes_\bfa i})).$$

\begin{lem}\label{homologyin} 
Suppose $(\bfa,\bfm)$ is a $j$-graded classical Rickard object in $\mathcal{T}$ such that $\bfa$ is Koszul with $\bfa^!$ concentrated in $k$-degree $0$ (that is, the pair $(\bfa,\bfa^!)$ satisfy the hypotheses in Setup \ref{gradconv}\eqref{conv2}). Set $\um = (\bfm,\bfm^{-1})$ and $X:=\bbH\bbT_\bfa(\um)$.
Then $\bbHH(\bfa,\bfm)$ is isomorphic to $\ffHH(\bbH\bbT_\bfa(\um))$ as  $ijk$-graded vector spaces, both being isomorphic to $$\bbH(\bigoplus_{s,t} (e_s \otimes e_s)(\bfa^{!} \otimes X)(e_t \otimes e_t))$$
with differential
$$a \otimes x \mapsto -\sum_{\rho \in \ttB^1} \left(a \rho^* \otimes \rho x -(-1)^{|x|} \rho^* a \otimes x \rho\right).$$
\end{lem}

\proof 
Since $\bbHH(\bfa,\bfm)\cong\bigoplus_{i\in \ZZ}  \HH(\bfa, \bfm^{\otimes_\bfa i})$ and each $\bfm^{\otimes_\bfa i}$ is a differential $jk$-graded $\bfa$-$\bfa$-bimodule, we can apply Theorem \ref{bimodulehochschild} to obtain an isomorphism
$$\bbHH(\bfa,\bfm) \cong \bigoplus_{i\in \ZZ} \bbH(\bigoplus_{s,t} e_s\bfa^!e_t \otimes e_t \bbH( \bfm^{\otimes_\bfa i})e_s) \cong \bbH(\bigoplus_{s,t} e_s\bfa^!e_t \otimes e_t \bbH( \bigoplus_{i\in \ZZ}\bfm^{\otimes_\bfa i})e_s) $$
where the differential on $\bigoplus_{s,t} e_s\bfa^!e_t \otimes e_t \bbH(\bigoplus_{i\in \ZZ} \bfm^{\otimes_\bfa i})e_s$ is as given in the statement of the lemma. Here we have used that homology and tensor products commute with direct sums.
Applying Theorem \ref{bimodulehochschild} to the $ijk$-graded $\bfa$-$\bfa$-bimodule (with trivial differential) $X$, we obtain the same result.
\endproof

Note that via this isomorphism, $\ffHH(\bbH\bbT_\bfa(\um))$ is equipped with a structure of associative algebra.
\begin{prop}\label{multstruc} Under the assumptions of Lemma \ref{homologyin}, the multiplicative structure on $\ffHH(\bbH\bbT_\bfa(\um))\cong\bbHH(\bfa,\bfm)$ is induced by the multiplicative structure on  $\bbH\bbT_\bfa(\um)$.
\end{prop}

\proof
By Proposition \ref{cup}, the multiplicative structure on $\bbHH(\bfa,\bfm)$ under the isomorphism to  $$\bbH(\bigoplus_{s,t} e_s\bfa^!e_t \otimes e_t \bbH( \bigoplus_{i\in \ZZ}\bfm^{\otimes_\bfa i})e_s)\cong \bbH(\bigoplus_{s,t} e_s\bfa^!e_t \otimes e_t  \bigoplus_{i\in \ZZ}\bbH(\bfm^{\otimes_\bfa i})e_s)$$ is given by $(a\otimes x)(b\otimes y) = ba \otimes xy$, where $xy$ is the multiplication induced in homology from the tensor product structure. At the same time, this is precisely the multiplication on $\bbH\bbT_\bfa(\um)$ and the claim follows.\endproof

\section{Representations of $\GL_2(F)$.} \label{GL2section}
Let $G=\GL_2(F)$. We study Ringel self-dual blocks of polynomial representations of $G$, where by a block of an abelian category $\cA$, we mean a Serre subcategory $\cB$ of $\cA$
minimal such that, given a pair of objects $L,M \in \cA$ with $\Ext^*_\cA(L,M) \neq 0$, the conditions $L \in \cB$ and $M \in \cB$ are equivalent.
According to \cite{EH} a block of polynomial representations of $G$ is  Ringel self-dual if and only if it has $p^l$ simple modules.

Let $\bfc$ be the finite dimensional algebra given by the quotient of the path algebra of $$\xymatrix{ \overset{1}{\bullet}
\ar@/^/[r]^{\eta} &\overset{2}{\bullet}\ar@/^/[l]^{\xi}
\ar@/^/[r]^{\eta} &\ar@/^/[l]^{\xi} \overset{3}{\bullet} &\cdots
&\overset{p-1}{\bullet}\ar@/^/[r]^{\eta} &
\ar@/^/[l]^{\xi} \overset{p}{\bullet},}$$ modulo the ideal
$$I=(\eta \xi e_p, \xi^2, \eta^2, \xi\eta +
\eta\xi).$$ 
The algebra $\bfc$ is $jk$-graded with $\eta$ and $\xi$ having $j$-degree $1$ and the whole algebra being concentrated in $k$-degree $0$.
It is a Ringel self-dual algebra with tilting bimodule $t$. Explicitly, $t$ can be defined as follows. We can realise the algebra $\bfc$ as an idempotent subquotient of the infinite-dimensional (non-unital) algebra  $Z$  given by the quiver
$$
\xymatrix{ \cdots & \overset{0}{\bullet} \ar@/^/[r]^{\eta}
& \ar@/^/[l]^{\xi} \overset{1}{\bullet} \ar@/^/[r]^{\eta}
&\overset{2}{\bullet}\ar@/^/[l]^{\xi} \ar@/^/[r]^{\eta}
&\ar@/^/[l]^{\xi} \overset{3}{\bullet} &\cdots
 },$$
modulo relations $\xi^2 = \eta^2 = \xi \eta + \eta \xi = 0$. Denote  by $\tau$ the algebra involution of $Z$ which sends vertex
$i$ to vertex $p-i$ and exchanges $\xi$ and $\eta$. Setting $$t = \sum_{1 \leq l \leq p, 0 \leq m \leq p-1} e_l Z e_m,$$
$t$ admits a natural left action of $\bfc$ by the subquotient $\bfc$ and a natural right action by twisting
the regular right $Z$-action by $\tau$. In this way, $t$ is naturally a $\bfc$-$\bfc$-bimodule.

We now let $\bft$ be a Rickard tilting complex representing $t$ for $\bfc$, and set $\bft^{-1} = \Hom_{\bfc}(\bft,\bfc)$ to be  its adjoint. It is then immediate that $(\bfc,\bft)$ is a classical Rickard object in $\mathcal{T}$. Indeed, $\bfc$ is an algebra, $\bft$ is projective on both sides and the natural morphism $\bfc\to \End_\bfc(\bft)$ is a quasi-isomorphism by Ringel self-duality.

 By \cite[Corollary 21]{MT2}, a block of polynomial representations of $G$ with $p^l$ simple modules is equivalent to the category of modules over the algebra $\bbO_{F,0} \bbO_{\bfc,t}^l(F,F)$. To compare to the notation used there, note that $\bbO_{F,0}$ simply picks out the algebra component of the resulting pair. By \cite[Lemma 30]{MT3}, there is an quasi-isomorphism $\bbO_{F,0} \bbO_{\bfc,\bft}^l(F,F) \to \bbO_{F,0} \bbO_{\bfc,t}^l(F,F)$ and quasi-isomorphic dg algebras share the same Hochschild cohomology, hence we define
$\hh_l$ to be the Hochschild cohomology of the algebra $\mathbb{O}_{F,0} \mathbb{O}_{\bfc,\bft}^l(F,F)$.
As we are ultimately interested in $\hh_l$ with $k$-grading given by the homological grading on Hochschild cohomology, 
we work with the gradings that suit this purpose, 
i.e. $\bfc$ is assumed to be concentrated in $k$-degree zero and $\bfc^!$ is assumed to be concentrated in positive $k$-degrees.

The aim of the rest of this article is to compute $\hh_l$.

\section{Reduction.} \label{computation}

The following Proposition demonstrates how our formalism of algebraic operators and homological duality 
reduce the computation of the algebra $\hh_l$ to the computation of the algebra 
$\ffHH(\bbH\bbT_{\bfc^!}(\ut^!))$, where $\ut^! = (\bft^!, \bft^{!-1})$ is the image of $(\bft, \bft^{-1})$ under Koszul duality.

\begin{prop} \label{stock}
We have $\hh_l \cong \fO_{F} \fO_{\ffHH(\bbH\bbT_{\bfc^!}(\ut^!))}^l (F[z,z^{-1}]).$
\end{prop}

\proof
We have algebra isomorphisms
\begin{equation*}\begin{split}
\hh_l  &\cong \bbHH \bbO_{F,0} \bbO_{\bfc,\bft}^l(F,F)  \\ 
&\cong \fO_F\bbHH( \bbO_{\bfc,\bft}^l(F,F))  \\ 
&\cong \fO_F \fO_{\bbHH(\bfc, \bft)}^l(\bbHH( F,F))\\
&\cong\fO_F \fO_{\bbHH(\bfc, \bft)}^l(F[z,z^{-1}])
\end{split}\end{equation*}
by Theorem \ref{iteration} and the observation that $\bbHH( F,0) \cong F$ and $\bbHH( F,F) \cong F[z,z^{-1}]$.

Rather than computing $\bbHH(\bfc, \bft)$ directly as $\bigoplus_i \HH(\bfc, \bft^i)$, 
we pull $\bfc$ through Koszul duality. We have derived equivalences (\cite[Appendix B]{MT4})
$$D(\bfc \lbigr) \cong D(\bfc^! \lbigr), \qquad D(\rbigr \bfc) \cong D(\rbigr \bfc^!)$$
$$\bfc \mapsto \bfc^{!*} \otimes_{\bfc^0} \bfc, \qquad \qquad \bfc \mapsto \bfc \otimes_{\bfc^0} \bfc^{!}.$$
Here $D(\bfc \lbigr)$ denotes the derived category of differential $jk$-bigraded left $\bfc$-modules and 
$D(\rbigr \bfc)$ denotes the derived category of differential $jk$-bigraded right $\bfc$-modules.
Putting these together (cf.~\cite[Theorem 2.1]{Rickard}) we have
\begin{equation*}\begin{split}
D(\bfc \bibigr \bfc ) &\cong  D( \bfc^! \bibigr \bfc^! )\\
\bfc &\mapsto \bfc^{!*} \otimes_{\bfc^0} \bfc \otimes_{\bfc^0} \bfc^!, \end{split}\end{equation*}
and since the equivalences $(- \otimes_{\bfc^!} \bfc^{!*} \otimes_{\bfc^0} \bfc , -\otimes_{\bfc} \bfc \otimes_{\bfc^0} \bfc^! )$ 
are adjoint equivalences (cf. \cite[Appendix B, \emph{Adjunction}]{MT4}) we have an isomorphism in the derived category between 
$\bfc^!$ and $\bfc^{!*} \otimes_{\bfc^0} \bfc \otimes_{\bfc^0} \bfc^!,$
Furthermore, by definition $\bft^!$ is the image of $\bft$ under the above equivalence. We thus have an isomorphism 
\begin{equation*}
\begin{split}
\bbHH(\bfc, \bft) & = \bigoplus_i \HH(\bfc, \bft^{\otimes_\bfc i})\\
&\cong \bigoplus_i \bbH\RHom_{\bfc \otimes \bfc^{\op}}(\bfc, \bft^{\otimes_\bfc  i})\\
&\cong \bigoplus_i \bbH\RHom_{\bfc^! \otimes \bfc^{! \op}}(\bfc^!,  \bfc^{!*} \otimes_{\bfc^0} \bft^{\otimes_\bfc i} \otimes_{\bfc^0} \bfc^!)\\
&\cong \bigoplus_i \bbH\RHom_{\bfc^! \otimes \bfc^{! \op}}(\bfc^!, (\bfc^{!*} \otimes_{\bfc^0} \bft \otimes_{\bfc^0} \bfc^!)^{\otimes_{\bfc^!}  i})\\
&\cong \bigoplus_i \bbH\RHom_{\bfc^! \otimes\bfc ^{! \op}}(\bfc^!, t^{!\otimes_{\bfc^!} i})\\
& = \bbHH(\bfc^!,\bft^!).
\end{split}
\end{equation*}
which implies  
$$\fO_F \fO_{\bbHH(\bfc,\bft)}^l(F[z,z^{-1}]) \cong \fO_F\fO_{\bbHH(\bfc^!,\bft^!)}^l(F[z,z^{-1}]).$$
But $\bbHH(\bfc^!,\bft^!)$ is isomorphic to $\ffHH(\bbH \bbT_{\bfc^!}(\ut^!))$ by Lemma \ref{homologyin},
which completes the proof of the Proposition.
\endproof

The above Proposition leaves us with the problem of computing $\ffHH(\bbH\bbT_{\bfc^!}(\ut^!))$ in the remaining sections.
We compute $\bbH\bbT_{\bfc^!}(\ut^!)$ in Section \ref{clubsuitsection}, 
then the Hochschild cohomology of the bimodules appearing in $\HH(\bfc^!, \bbH(\bft^{!\otimes_{\bfc^!} i}))$ for various $i$ in Section \ref{explicitbimod},  
and finally infer the multiplication on $\ffHH( \bbH\bbT_{\bfc^!}(\ut^!))$ 
from that on $\bbH\bbT_{\bfc^!}(\ut^!)$ in Section \ref{spadesuitsection}.

\section{The algebra $\HTct$.} \label{clubsuitsection}

In this section we compute the algebra structure of $\HTct:= \mathbb{H}\bbT_{\bfc^!}(\ut^!)$, which entwines the algebra $\bfc^!$, its dual, its tilting bimodule,
and a preprojective algebra $\Theta$ in a subtle way. We do this by first computing $\HTct^-:= \bbH(\bbT_{\bfc^!}(\bft^{!-1}))$ and $\HTct^+:= \bbH(\bbT_{\bfc^!}(\bft^{!}))$ separately and then investigating their interaction.

\subsection{The algebras $\Omega$ and $\Theta$.}
We first need some notation.
The algebra $\bfc$ has generators $\xi$ and $\eta$, and its Koszul dual $\bfc^!=\Omega$ has dual generators $x$ and $y$;
The quiver of $\Omega$ is given by
$$
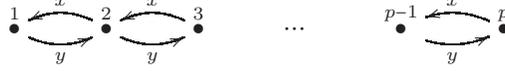
\xymatrix{ \overset{1}{\bullet} \ar@/_/[r]_{y}
&\overset{2}{\bullet}\ar@/_/[l]_{x} \ar@/_/[r]_{y}
&\ar@/_/[l]_{x} \overset{3}{\bullet} &\cdots & \overset{p-1}{\bullet} \ar@/_/[r]_{y}
&\ar@/_/[l]_{x} \overset{p}{\bullet}
 }$$
and the relations for $\Omega$ are $xye_1=0$ and $xy=yx$.
Since $\xi$ and $\eta$ were assumed to have $j$-degree $1$ and $k$-degree $0$, $x$ and $y$ now both have $j$-degree $-1$ and $k$-degree $1$.
For notational convenience we use a different convention for the direction of arrows in $\Omega$ than we used in our previous article \cite{MT3}. We denote by $e_i$ the idempotent corresponding to vertex $i$.

Note that morphisms from $\Omega e_i$ to $\Omega e_l$ are of the form $\cdot x^{l-i}(xy)^s$ if $i\leq l$, where $0\leq s \leq i-1$, or of the form
$\cdot y^{i-l}(xy)^s$ if $i\geq l$, where $0\leq s \leq l-1$. Such morphisms have $(j,k)$ degree $(-(l-i+2s),l-i+2s)$ and $(-(i-l+2s),i-l+2s)$ respectively.

The algebra $\Omega$ has a simple preserving duality, interchanging $x$ and $y$. It is quasi-hereditary (with uniserial standard modules $\Delta_i=\Omega e_i/\Omega e_{i-1}\Omega e_i$) and Ringel self-dual (by \cite[Theorem 1, Example 19]{Maz} and Ringel self-duality of $\bfc$), and its tilting (bi-)module is easily seen to be isomorphic to $\Omega e_p \Omega$.
This bimodule is self-dual via the isomorphism \begin{equation}\label{omegaiso}\Omega e_p\Omega\cong (\Omega e_p\Omega)^*\langle 2-2p\rangle [2-2p]\end{equation} induced by the 
symmetric associative nondegenerate bilinear form 
$$\Omega e_p \Omega \otimes \Omega e_p \Omega \twoheadrightarrow F,$$
sending $e_s x^d y^e e_t \otimes e_{s'} x^{d'} y^{e'} e_{t'}$ to $1$ if $s=t'$, $t=s'$, and $d+d'=e+e'=p-1$, and to zero otherwise. The degree shift comes from the bimodule socle of $\Omega e_p\Omega$ (which is given by $e_p y^{p-1}x^{p-1}e_p$) having $(j,k)$-degree $(2-2p,2p-2)$, and thus the bimodule top of  $(\Omega e_p\Omega)^*$ having $(j,k)$-degree $(2p-2,2-2p)$. Thus, with our grading conventions from Section \ref{grad}, $(\Omega e_p\Omega)^*\langle 2-2p\rangle [2-2p]$ indeed has top in degree $(0,0)$.
We furthermore claim that $(\Omega, \Omega e_p\Omega)$ is a Rickard object. Indeed, $\Omega$ is an  algebra, $\Omega e_p\Omega$ is perfect both as a left and as a right $\Omega$-module, and $\RHom_\Omega(\Omega e_p\Omega,\Omega e_p\Omega)$ is in fact isomorphic to $\Omega$, since $\Ext^i_\Omega(\Omega e_p\Omega,\Omega e_p\Omega)=0$ for $i>0$ thanks to $\Omega e_p\Omega$ being the tilting module for a Ringel self-dual quasi-hereditary algebra.

We define the algebra $\Theta$ to be the quotient $\Omega/\Omega e_p \Omega$, where $e_i$ denotes the idempotent at vertex $i$.
The algebra $\Theta$ is called the \emph{preprojective algebra of type $A_{p-1}$}. Let $\sigma$ be the involution of $\Theta$ which switching $e_s$ and $e_{p-s}$, and $x$ and $y$.
Then $\Theta$ is a self-injective algebra with Nakayama automorphism $\sigma$. Indeed we have an isomorphism of $\Theta$-$\Theta$-bimodules
\begin{equation}\label{thetaiso}\Theta^\sigma \to \Theta^*\langle 2-p\rangle [2-p]\colon \; e_s \mapsto e_s(y^{p-s-1}x^{s-1})^*e_{p-s}.\end{equation}
Indeed it is easy to check  that this is an isomorphism of ungraded bimodules, and the degree $(0,0)$-part of $\Theta^*\langle 2-p\rangle [2-p]$ is, according to our grading conventions from Section \ref{grad}, equal to $(\Theta^*)^{p-2,2-p}=(\Theta^{2-p,p-2})^*$, which indeed contains the element $e_s(y^{p-s-1}x^{s-1})^*e_{p-s}$.

Viewed as a tilting complex of ungraded left $\Omega$-modules, $\bft^{!-1}$ is quasi-isomorphic to the direct sum of 
$$\Omega e_p \to 0  $$ with the direct sum over $l=1, \ldots, p$
of two term complexes
$$\Omega e_p \overset{\cdot y^{l}}{\to} \Omega e_{p-l}$$
by \cite[Lemma 34]{MT3}.
By \cite[Lemma 37(iv)]{MT3}, the right action of $\Omega$ on these complexes is given by the action of $e_{l}xe_{l+1}$ respectively $e_{l}ye_{l-1}$  (whenever none of involved idempotents are $e_p$) as

\begin{equation}\label{y1}\xymatrix{
\Omega e_p  \ar^{\cdot y^{l}}[r] \ar_1[d] & \Omega e_{p-l} \ar^{\cdot y}[d] &&&\Omega e_p\ar^{\cdot y^{l}}[r] \ar_{\cdot xy}[d] & \Omega e_{p-l} \ar^{\cdot x}[d]  \\
\Omega e_p \ar^{\cdot y^{l+1}}[r] & \Omega e_{p-l-1} &&&\Omega e_p \ar^{\cdot y^{l-1}}[r] & \Omega e_{p-l+1} ,
 }\end{equation}
while the action of the elements $e_{p-1}xe_p$ and $e_pye_{p-1}$  is given by
\begin{equation}\label{x1}\xymatrix{
\Omega e_{p} \ar^{\cdot y^{p-1}}[r] \ar[d]_1 & \Omega e_{1} &&& \Omega e_p\ar_{\cdot xy}[d]& \\
\Omega e_p & &&&\Omega e_p  \ar^{\cdot y^{p-1}}[r] &\Omega e_{1}}\end{equation}
respectively.

Taking the adjoint and applying our simple-preserving duality, $\bft^!$ is, as a tilting complex of ungraded left $\Omega$-modules, quasi-isomorphic to the direct sum of 
$$0\to \Omega e_p  $$ with the direct sum over $l=1, \ldots, p$
of two term complexes
$$\Omega e_{p-l} \overset{\cdot x^{l}}{\to} \Omega e_p$$
with the right action of the generators  $e_{l}xe_{l+1}$ respectively $e_{l}ye_{l-1}$ (whenever none of involved idempotents are $e_p$) given by 
\begin{equation}\label{eq11}\xymatrix@C=15pt{
 \Omega e_{p-l}\ar^{\cdot y}[d]\ar^{\cdot x^l}[r] & \Omega e_p \ar^{\cdot xy}[d] &&&&
 \Omega e_{p-l}\ar^{\cdot x}[d]\ar^{\cdot x^l}[r] &\Omega e_p \ar^{\cdot 1}[d] \\
\Omega e_{p-l-1}\ar^{\cdot x^{l+1}}[r] & \Omega e_p  &&&&
 \ar^{\cdot x^{l-1}}[r]\Omega e_{p-l+1}& \Omega e_p 
}
\end{equation}
respectively, while the action of the elements $e_{p-1}xe_p$ and $e_pye_{p-1}$  is given by
\begin{equation}\label{eq11p}\xymatrix@C=15pt{
\ar^{\cdot x^{p-1}}[r]\Omega e_{1}& \Omega e_p \ar^{xy}[d] &&&&&
\Omega e_p \ar^{\cdot 1}[d]   \\
&\Omega e_p & & &  &
\Omega e_{1}\ar^{\cdot x^{p-1}}[r]&\Omega e_p  
}
\end{equation}
respectively.

To ease notation, we will, in the remainder of the article, write $\bft^{!i}$ for $(\bft^{!})^{\otimes_{\bfc^!} i}=(\bft^{!})^{\otimes_{\Omega} i}$ and $\bft^{!-i}$ for $(\bft^{!-1})^{\otimes_{\bfc^!} i}=(\bft^{!-1})^{\otimes_\Omega i}$ for $i>0$.

\subsection{Recollections of the homology $\bbH(\bbT_{\bfc^!}(\bft^{! -1}$)).}\label{hctinvrecall}

Recall that, given a collection $(M_s,f_s)_s$ where the $M_s$ are differential $(j,k)$-bigraded $\Omega$-modules, and the $f_s$ are morphisms of differential $(j,k)$-bigraded $\Omega$-modules (preserving both $j$- and $k$-degrees), such that the sequence $$\cdots M_s\overset{f_s}{\to}M_{s+1} \overset{f_{s+1}}{\to}M_{s+1}\cdots$$
is a complex of $(j,k)$-bigraded vector spaces, we can associate a 
differential $(j,k)$-bigraded $\Omega$-module, namely the iterated cone of the family of morphisms $(f_s)_s$. In particular, if 
$$f : M_{-s} \overset{f_{-s}}{\to} M_{-s+1} \overset{f_{-s+1}}{\to}\cdots \overset{f_{-1}}{\to} M_0$$
is a complex of $(j,k)$-bigraded vector spaces, the homology of this sequence is given by $\bigoplus_{i=0}^s H^i(f)[i]$.

We now summarise the results of \cite[Section 8]{MT3}, recalling that there $x$ and $y$ were interchanged, and given $j$-degree $1$, so in particular, all shifts in $j$-degree from \cite{MT3} appear as the negative here.

Consider the family $(f^l)_{l=1,\ldots, p}$ of morphisms differential $(j,k)$-bigraded $\Omega$-modules 
$$f_l\colon \Omega e_p\langle -l \rangle [-l] \to \Omega e_{p-l} $$
given by right multiplication with $y^{l}$ for $l=1,\ldots,p-1$ and by the zero map 
$$f_p\colon \Omega e_p\langle -p \rangle [-p] \to 0. $$
By \cite[Lemma 34]{MT3}, the differential $(j,k)$-bigraded $\Omega$-module $\bft^{!-1}$  is quasi-isomorphic to the cone $X^{-1} $ of the direct sum $\bigoplus_{l=1}^pf_l$ of these morphisms. By \cite[Lemma 37 (iv)]{MT3}, $X^{-1} $ has homology isomorphic to $\Omega \langle -p \rangle[1-p] \oplus\Theta^\sigma$, coming from a direct sum over $l$ of  exact sequences of $j$-graded $\Omega$-modules \cite[Lemma 35]{MT3}
$$0\to\Omega e_l\langle -p \rangle\to \Omega e_p\langle -l \rangle \to \Omega e_{p-l}\to \Theta e_{p-l} \to 0,$$
to which the homology of the isolated summand coming from $f_p$ is added. The right  action of $\Omega$ is induced by the diagrams \eqref{y1} and \eqref{x1}.

For $i>1$, by \cite[Lemma 38]{MT3},$\bft^{! -i}$ is quasi-isomorphic to the direct sum $X^{-i}=\bigoplus_{l=1}^pX^{-i}e_l$, where $X^{-i}e_l$ is the iterated cone of the $i+1$-term sequence $g^l$ of morphisms of differential $(j,k)$-bigraded $\Omega$-modules (each sequence being  a complex of $(j,k)$-bigraded vector spaces), where $g^p$ is given by 
$$\Omega e_p \langle -ip \rangle[-ip]  \to 0\to\cdots \to 0,$$
 and, for $l=1,\ldots ,p-1$,  $g^l$ is the sequence
$$ \cdots\overset{\cdot (xy)^{l}}{\rightarrow} \Omega e_p \langle l-3p \rangle [l-3p]\overset{\cdot (xy)^{p-l}}{\rightarrow}\Omega e_p \langle -(l+p) \rangle [-(l+p)] \overset{\cdot (xy)^{l}}{\rightarrow} \Omega e_p\langle l-p \rangle [l-p]\overset{\cdot y^{p-l}}{\rightarrow} \Omega e_l$$
if $i$ is even and
$$ \cdots\overset{\cdot (xy)^{p-l}}{\rightarrow} \Omega e_p \langle -(l+2p) \rangle [-(l+2p)]\overset{\cdot (xy)^{l}}{\rightarrow}\Omega e_p \langle l-2p \rangle [l-2p] \overset{\cdot (xy)^{p-l}}{\rightarrow} \Omega e_p\langle -l \rangle [-l]\overset{\cdot y^{l}}{\rightarrow} \Omega e_{p-l}$$
if $i$ is odd.

Furthermore $X^{-i}$ has homology
\begin{equation*}
\begin{split}
\bbH(\bft^{! -i})&\cong\Omega \langle -ip \rangle[i(1-p)] \oplus \Theta^{\sigma} \langle -(i-1)p \rangle[(i-1)(1-p)]
\oplus ...\oplus \Theta^{\sigma^i} \langle 0 \rangle[0] \\
&\cong \Omega \langle -ip \rangle[i(1-p)] \oplus \bigoplus_{j=1}^i  \Theta^{\sigma^j} \langle -(i-j)p \rangle[(i-j)(1-p)]. \end{split}\end{equation*}

The structure of $\HTct^-= \bbH(\bbT_{\bfc^!}(\bft^{!-1}))$ as a $k$-graded $\Omega$-$\Omega$-bimodule is therefore given by
$$
\xymatrix@C=10pt@R=3pt{
&&&&& \Omega  [0] & & &      \\
&&&& \Omega [1-p]  & \Theta^\sigma   [0]   &&&&&  \\
&&& \Omega [2-2p] & \Theta^\sigma [1-p] & \Theta [0]   &&&&&   \\
&& \Omega [3-3p] & \Theta^\sigma [2-2p] & \Theta [1-p] & \Theta^\sigma [0]  &&&&&     \\
& \Omega [4-4p] & \Theta^\sigma [3-3p] & \Theta [2-2p] & \Theta^\sigma [1-p] & \Theta [0] &&&&&      \\
& & & ... & &
}$$
and the structure of $\HTct^-$ as a $j$-graded $\Omega$-$\Omega$-bimodule is given by
$$
\xymatrix@C=10pt@R=3pt{
&&&&& \Omega  \langle 0 \rangle & & &      \\
&&&& \Omega \langle -p \rangle  & \Theta^\sigma    \langle 0 \rangle  &&&&&  \\
&&& \Omega \langle -2p \rangle & \Theta^\sigma \langle -p \rangle & \Theta  \langle 0 \rangle  &&&&&   \\
&& \Omega  \langle -3p \rangle & \Theta^\sigma \langle -2p \rangle & \Theta \langle -p \rangle & \Theta^\sigma  \langle 0 \rangle &&&&&     \\
& \Omega  \langle -4p \rangle & \Theta^\sigma \langle -3p \rangle & \Theta \langle -2p \rangle & \Theta^\sigma \langle -p \rangle & \Theta \langle 0 \rangle. &&&&&      \\
& & & ... & &
}$$

By \cite[Theorem 32]{MT3}, $\HTct^-$ is is isomorphic to the tensor algebra $\bbT_\Omega(\Theta^\sigma) \otimes F[\xi]$ where $\xi$ is a variable of $j$-degree $-p$ and $k$-degree $p-1$, so that 
$\Omega \xi \cong \Omega\langle -p \rangle [1-p]$.

\subsection{Homology of the bimodules $\bft^{! i}$ for $i>0$.}
By definition, $\bft^!=\Hom_\Omega(\bft^{! -1}, \Omega)$, so using our simple preserving duality,  $\bft^!$ is quasi-isomorphic to the direct sum $X^1=\bigoplus_{l=1}^pX^1e_l$, where $X^1e_l[1]$ is the cone of the morphism 
$$\Omega e_{p-l}\overset{\cdot x^l}{\to}\Omega e_p\langle l \rangle [l] $$
for $l=1,\ldots,p-1$, and of
$$0\to\Omega e_p\langle p \rangle [p] $$
for $l=p$.
The homology of $X^1e_l$ is easily seen to be $\Omega e_p\Omega e_l\langle p \rangle [p-1]$, so $X^1$ is in fact quasi-isomorphic $\Omega e_p\Omega \langle p \rangle [p-1]$.

For $\bft^{! 2}$, we similarly see that this is quasi-isomorphic to 
$X^2=\bigoplus_{l=1}^pX^2e_l$ where $X^2e_l[2]$ is the iterated cone of the sequence of morphisms
$$ \Omega e_l\overset{\cdot x^{p-l}}{\to}  \Omega e_p \langle p-l \rangle [p-l] \overset{\cdot (xy)^l}{\to} \Omega e_p \langle p+l \rangle [p+l]$$
for $l=1,\ldots,p-1$, and of
$$0\to0\to \Omega e_p\langle 2p \rangle [2p] $$
for $l=p$.

Using that $\Omega e_p\Omega$ is the tilting module for the Ringel self-dual algebra $\Omega$ and the isomorphism given in \eqref{omegaiso}, we have a sequence of  ungraded $\Omega$-$\Omega$-bimodule isomorphisms,
\begin{equation}\label{tiltingsquare}
\begin{split}
\Omega e_p\Omega\otimes_\Omega \Omega e_p\Omega &\cong (\Omega e_p\Omega)^*\otimes_\Omega \Omega e_p\Omega \\
&\cong \Hom_F( \Hom_F((\Omega e_p\Omega)^*\otimes_\Omega \Omega e_p\Omega ,F),F)\\
&\cong\Hom_F(\Hom_\Omega(\Omega e_p\Omega, \Omega e_p\Omega),F)
\\&\cong \Omega^*
\end{split}
\end{equation}
 as an ungraded $\Omega$-$\Omega$-bimodule. Explicitly, denoting by $\langle - , - \rangle$ the pairing obtained from \eqref{omegaiso}, an isomorphism is given by the assignment $$u\otimes v\mapsto \left(w \mapsto \langle u,vw\rangle\right). $$ Thus we already know that the homology of $X^2$ is isomorphic to $\Omega^*$ as an ungraded $\Omega$-$\Omega$-bimodule, and we only need to determine the gradings. Direct computation shows that the sequence of  morphisms 
$$ \Omega e_l\overset{\cdot x^{p-l}}{\to}  \Omega e_p \langle p-l \rangle [p-l] \overset{\cdot (xy)^l}{\to} \Omega e_p \langle p+l \rangle [p+l]$$
 indeed has homology $\Omega^*e_l\langle 2\rangle[2]$ in the last place, via the isomorphism in homology induced by the morphism of left $\Omega$-modules $\Omega e_p \to \Omega^*e_l$ which sends $e_p$ to $(e_ly^{l-1}x^{p-1}e_p)^*$. Hence the homology of $X^2e_l$ is given by $\Omega^*e_l\langle 2\rangle[0]$.

Again using the simple preserving duality, we see that for $i>2$, $t^{!i}$ is quasi-isomorphic to 
$X^{i}=\bigoplus_{l=1}^pX^{i}e_l$, where $X^{i}e_l[i]$ is the iterated cone of the $i+1$-term sequence $\tilde g^l$ of morphisms of differential $(j,k)$-bigraded $\Omega$-modules (each sequence being  a complex of $(j,k)$-bigraded vector spaces), where $\tilde g^p$ is given by 
$$  0\to\cdots \to 0\to \Omega e_p \langle ip \rangle[ip],$$
 and, for $l=1,\ldots ,p-1$,  $\tilde g^l$ is the sequence
$$ 
\Omega e_l\overset{\cdot x^{p-l}}{\rightarrow}  \Omega e_p\langle p-l \rangle[p-l]\overset{\cdot (xy)^{l}}{\rightarrow}\Omega e_p \langle l+p \rangle [l+p] \overset{\cdot (xy)^{p-l}}{\rightarrow}\Omega e_p \langle 3p-l \rangle [3p-l]\overset{\cdot (xy)^{l}}{\rightarrow} \cdots
$$
if $i$ is even and
$$  \Omega e_{p-l}\overset{\cdot x^{l}}{\rightarrow}\Omega e_p\langle l \rangle [l]\overset{\cdot (xy)^{p-l}}{\rightarrow} \Omega e_p \langle 2p-l \rangle [2p-l] \overset{\cdot (xy)^{l}}{\rightarrow}\Omega e_p \langle l+2p \rangle [l+2p]\overset{\cdot (xy)^{p-l}}{\rightarrow} \cdots$$
if $i$ is odd.

Since $\Omega$ is a finite dimensional algebra which has finite global dimension (as the Koszul dual of a finite-dimensional algebra, alternatively, as a quasi-hereditary algebra),  $\Omega^* \otimes^{\mathbb{L}}_\Omega -$ is a Serre functor on $D^b(\Omega)$ by \cite[4.6]{Ha}. 
Hence we have, in the ungraded setting,  a quasi-isomorphism 
$$\bft^{! -i} =\Hom_{\Omega}(\bft^{! i}, \Omega)\overset{qim}{\rightarrow} \Hom_{\Omega}(\Omega, \Omega^* \otimes_\Omega \bft^{! i})^*$$ 
where we have used that 
$$\Hom_{\Omega}(\Omega, \Omega^* \otimes^{\mathbb{L}}_\Omega \bft^{! i})^*\cong\Hom_{\Omega}(\Omega, \Omega^* \otimes_\Omega \bft^{! i})^*$$
 since $\bft^{! i}$ is projective as a left $\Omega$-module. Thus $\bft^{! -i}$ is quasi-isomorphic to $$(\Omega^* \otimes_\Omega \bft^{! i})^* = (\bft^{!  i+2})^*.$$

Putting in gradings, this gives a quasi-isomorphism between 
$$\bft^{! -i}  \overset{qim}{\rightarrow}\Hom_{\Omega}(\Omega, \Omega^* \otimes_\Omega \bft^{! i})^*\cong \Hom_{\Omega}(\Omega, \bft^{!  i+2}\langle -2 \rangle [0] )^* \cong (\bft^{!  i+2 })^*\langle 2 \rangle [0] ,$$ or, equivalently,  for $i\geq 2$, a quasi-isomorphism
$$\bft^{! i} \overset{qim}{\rightarrow} (\bft^{! -(i-2)}\langle -2\rangle [0])^* = (\bft^{! -(i-2)})^*\langle 2\rangle [0].$$

Therefore,
\begin{equation*}
\begin{split}
 \bbH&(\bft^{! i}) \cong (\bbH(\bft^{! -(i-2)}))^*\langle 2\rangle \\&\cong ( \Omega \langle -(i-2)p\rangle[(i-2)(1-p)]  \oplus \Theta^\sigma \langle -(i-3)p\rangle[(i-3)(1-p)]  \oplus \cdots \oplus \Theta^{\sigma^{i-2}} \langle 0 \rangle [0])^* \langle 2 \rangle\\
& \cong (\Omega^*\langle (i-2)p\rangle[(i-2)(p-1)]\oplus \Theta^{\sigma *} \langle (i-3)p\rangle[(i-3)(p-1)]  \oplus \cdots \oplus \Theta^{\sigma^{i-2}*} \langle 0 \rangle [0])\langle 2 \rangle\\
&\cong \Omega^*\langle 2+(i-2)p\rangle[(i-2)(p-1)]\oplus \Theta^{\sigma *} \langle 2+(i-3)p\rangle[(i-3)(p-1)]  \oplus \cdots \oplus \Theta^{\sigma^{i-2}*} \langle 2\rangle [0]\\
& \cong \Omega^*\langle 2+(i-2)p\rangle[(i-2)(p-1)]\oplus 
\bigoplus_{j=1}^{i-2}\Theta^{\sigma^{j} *}\langle 2+(i-2-j)p\rangle[(i-2-j)(p-1)]. 
\end{split}
\end{equation*}
Using $\Theta^{*} \cong \Theta^\sigma \langle p-2 \rangle [p-2]$ coming from the isomorphism \eqref{thetaiso} and the fact that the involution $\sigma$ of $\Theta$ induces an isomorphism of bimodules ${}^\sigma\Theta \cong \Theta^{\sigma^{-1}}\cong \Theta^\sigma$, we obtain
\begin{equation}\label{plushomiso}
\begin{split}
\bbH(\bft^{! i})&
\cong 
  \Omega^*\langle 2+(i-2)p\rangle[(i-2)(p-1)]\oplus 
  \bigoplus_{j=1}^{i-2}\Theta^{\sigma^{j+1}}\langle (i-1-j)p\rangle[(i-1-j)(p-1)-1]\\
  &\cong  \Omega^*\langle 2+(i-2)p\rangle[(i-2)(p-1)]\oplus   \Theta^{} \langle (i-2)p\rangle[(i-2)(p-1)-1]  \oplus \cdots\\& \qquad\cdots \oplus \Theta^{\sigma^{i-1}} \langle p\rangle [p-2].
\end{split}
\end{equation}

Explicitly, the generator $e_p$ of the rightmost copy of $\Omega e_p$ in $X^ie_l$ corresponds to the element $(e_lx^{p-1} y^{l-1}e_p)^* \in \Omega^*$ in homology. The homology class of an element $ue_p$ in a middle term of the form $\Omega e_p$ in $X_i$ annihilated by the morphism given by right multiplication by $(xy)^l$ corresponds to the element $u'e_l$ in $\Theta$ such that $u'x^{p-l}e_p=ue_p$.

Hence the structure of $\HTct^+$ as a $k$-graded $\Omega$-$\Omega$-bimodule is given by
$$
\xymatrix@C=10pt@R=3pt{
&&&& & ... & & &\\
&&&&       &\Theta [p-2]      & \Theta^\sigma [2p-3] &  \Theta [3p-4] & \Omega^* [3p-3] \\
&&&&            & \Theta^\sigma [p-2] &  \Theta [2p-3] & \Omega^* [2p-2] &\\
&&&& & \Theta [p-2] & \Omega^* [p-1] && \\
&&&& & \Omega^*  [0] & && \\
&&&& \Omega e_p \Omega [p-1] & & &\\
&&& \Omega   & & &   &   \\
}$$
while the structure of $\HTct^+$ as a $j$-graded $\Omega$-$\Omega$-bimodule is given by
$$
\xymatrix@C=10pt@R=3pt{
&&&& & ... & & \\
&&&&        &\Theta \langle p \rangle  & \Theta^\sigma \langle 2p\rangle &  \Theta     \langle 3p\rangle      & \Omega^* \langle 2+3p \rangle\\
&&&&         & \Theta^\sigma \langle p\rangle &  \Theta     \langle 2p\rangle      & \Omega^* \langle 2+2p \rangle &\\
&&&&  & \Theta \langle p \rangle & \Omega^* \langle 2+p \rangle & &\\
&&&&  & \Omega^*  \langle 2\rangle & & &\\
&&&& \Omega e_p \Omega \langle p \rangle & & &\\
&&& \Omega   & & &    &.  \\
}$$

\subsection{The product on $\HTct$.}
We now investigate the algebra structure on $\HTct$. We recall from  \cite[Theorem 32]{MT3} (or Section \ref{hctinvrecall})
that $\HTct^-$ is nothing but the tensor algebra $\bbT_\Omega(\Theta^\sigma) \otimes F[\xi]$ for  a variable $\xi$ of $j$-degree $-p$ and $k$-degree $p-1$. 

In order to determine the products of two elements in $\HTct^+$, or mixed products between $\HTct^+$ and $\HTct^-$, we use an explicit right $\Omega$-module structure on the one-sided tilting comlpexes described in the previous section.

In \cite[Lemma 38 (ii), equations (9) and (10)]{MT3}, we gave the description of the right $\Omega$-structure on $X^{ -i}$ in the example of $i$ odd and not involving the $p$th summand. For completeness, we include a full description here.  The action of the generators  $e_{l}xe_{l+1}$ respectively $e_{l}ye_{l-1}$ is induced by the diagrams
\begin{equation}\label{eq1}\xymatrix@C=15pt{
\Omega e_p \ar^{\cdot 1}[d]\ar^{\cdot(xy)^{l}}[r]& \cdots \ar^{\cdot(xy)^{l}}[r]& \Omega e_p  \ar^{\cdot xy}[d] \ar^{\cdot(xy)^{p-l}}[r] &  \Omega e_p \ar^{\cdot 1}[d] \ar^{\cdot y^l}[r]&\Omega e_{p-l}\ar^{\cdot y}[d] &
\Omega e_p  \ar^{\cdot(xy)^{l}}[r]\ar^{\cdot xy}[d]& \cdots\ar^{\cdot(xy)^{l}}[r] & \Omega e_p  \ar^{\cdot 1}[d] \ar^{\cdot(xy)^{p-l}}[r] &  \Omega e_p \ar^{\cdot xy}[d] \ar^{\cdot y^l}[r]&\Omega e_{p-l}\ar^{\cdot x}[d] \\
\Omega e_p \ar^{\cdot(xy)^{l+1}}[r]&\cdots \ar^{\cdot(xy)^{l+1}}[r]&  \Omega e_p \ar^{\cdot(xy)^{p-l-1}}[r]&  \Omega e_p  \ar^{\cdot y^{l+1}}[r]&\Omega e_{p-l-1}&
\Omega e_p \ar^{\cdot(xy)^{l-1}}[r] &\cdots \ar^{\cdot(xy)^{l-1}}[r]&  \Omega e_p \ar^{\cdot(xy)^{p-l+1}}[r]&  \Omega e_p  \ar^{\cdot y^{l-1}}[r]&\Omega e_{p-l+1}
}
\end{equation}
for $i$ odd and
\begin{equation}\label{eq2}\xymatrix@C=15pt{
\Omega e_p  \ar^{\cdot(xy)^{l}}[r]\ar^{\cdot 1}[d]& \cdots\ar^{\cdot(xy)^{p-l}}[r] & \Omega e_p \ar^{\cdot(xy)^{l}}[r] \ar^{\cdot 1}[d] &  \Omega e_p\ar^{\cdot xy}[d] \ar^{\cdot y^{p-l}}[r]&\Omega e_{l}\ar^{\cdot x}[d]&
\Omega e_p  \ar^{\cdot(xy)^{l}}[r]\ar^{\cdot xy}[d]& \cdots \ar^{\cdot(xy)^{p-l}}[r]& \Omega e_p \ar^{\cdot(xy)^{l}}[r] \ar^{\cdot xy}[d] &  \Omega e_p\ar^{\cdot 1}[d] \ar^{\cdot y^{p-l}}[r]&\Omega e_{l}\ar^{\cdot y}[d]\\
\Omega e_p  \ar^{\cdot(xy)^{l+1}}[r] &\cdots\ar^{\cdot(xy)^{p-l-1}}[r]& \Omega e_p \ar^{\cdot(xy)^{l+1}}[r]&  \Omega e_p  \ar^{\cdot y^{p-l-1}}[r]&\Omega e_{l+1}&
\Omega e_p  \ar^{\cdot(xy)^{l-1}}[r] &\cdots\ar^{\cdot(xy)^{p-l+1}}[r]& \Omega e_p\ar^{\cdot(xy)^{l-1}}[r]&  \Omega e_p  \ar^{\cdot y^{p-l+1}}[r]&\Omega e_{l-1}
}
\end{equation}
for $i$ even, wherever this makes sense for $l$ (i.e. the larger value being less than or equal to $p-1$). The action of the elements $e_{p-1}xe_p$ and $e_pye_{p-1}$ is induced by
\begin{equation}\label{eq1p}\xymatrix@C=15pt{
\Omega e_p \ar^{\cdot 1}[d]\ar^{\cdot(xy)}[r]& \cdots \ar^{\cdot(xy)^{p-1}}[r]& \Omega e_p   \ar^{\cdot(xy)}[r] &  \Omega e_p \ar^{\cdot y^{p-1}}[r]&\Omega e_{1} &
\Omega e_p \ar^{\cdot xy}[d]&  &  & & \\
\Omega e_p &&  &  &&
\Omega e_p \ar^{\cdot(xy)^{p-1}}[r] &\cdots \ar^{\cdot(xy)^{p-1}}[r]&  \Omega e_p \ar^{\cdot(xy)}[r]&  \Omega e_p  \ar^{\cdot y^{p-1}}[r]&\Omega e_{1}
}
\end{equation}
for $i$ odd and 
\begin{equation}\label{eq2p}\xymatrix@C=15pt{
\Omega e_p  \ar^{\cdot(xy)^{p-1}}[r]\ar^{\cdot 1}[d]& \cdots\ar^{\cdot(xy)}[r] & \Omega e_p \ar^{\cdot(xy)^{p-1}}[r] &  \Omega e_p\ar^{\cdot y}[r]&\Omega e_{p-1}&
\Omega e_p  \ar^{\cdot xy}[d]&& & &\\
\Omega e_p && &  &&
\Omega e_p  \ar^{\cdot(xy)^{p-1}}[r] &\cdots\ar^{\cdot(xy)}[r]& \Omega e_p\ar^{\cdot(xy)^{p-1}}[r]&  \Omega e_p  \ar^{\cdot y}[r]&\Omega e_{p-1}
}
\end{equation}
for $i$ even.

Similarly, the right $\Omega$-structure on $\bft^{! i}$ is generated by the action of $e_{l}xe_{l+1}$ respectively $e_{l}ye_{l-1}$ induced from the morphism of complexes
\begin{equation}\label{eq3}\xymatrix@C=15pt{
\Omega e_{p-l} \ar^{\cdot x^l}[r]\ar^{\cdot y}[d]& \Omega e_p \ar^{\cdot (xy)^{p-l}}[r]\ar^{\cdot xy}[d]   &   \Omega e_p \ar^{\cdot 1}[d] \ar^{\cdot (xy)^{l}}[r]& \ar^{\cdot (xy)^{l}}[r]\cdots&\Omega e_{p}\ar^{\cdot xy}[d]&
\Omega e_{p-l} \ar^{\cdot x^l}[r]\ar^{\cdot x}[d]& \Omega e_p \ar^{\cdot (xy)^{p-l}}[r]\ar^{\cdot 1}[d]  &   \Omega e_p \ar^{\cdot xy}[d]\ar^{\cdot (xy)^{l}}[r] & \cdots\ar^{\cdot(xy)^{l} }[r]&\Omega e_{p}\ar^{\cdot 1}[d]\\
\Omega e_{p-l-1}  \ar^{\cdot x^{l+1}}[r]& \Omega e_p\ar^{\cdot (xy)^{p-l-1}}[r] &  \Omega e_p\ar^{\cdot (xy)^{l+1}}[r]&\cdots   \ar^{\cdot (xy)^{l+1}}[r]& \Omega e_{p}&
\Omega e_{p-l+1}  \ar^{\cdot x^{l-1}}[r]& \Omega e_p \ar^{\cdot (xy)^{p-l+1}}[r]&  \Omega e_p \ar^{\cdot (xy)^{l-1}}[r]&\cdots \ar^{\cdot (xy)^{l-1}}[r]  & \Omega e_{p}
}
\end{equation}
for $i$ odd and 
\begin{equation}\label{eq4}\xymatrix@C=15pt{
\Omega e_{l}  \ar^{\cdot x^{p-l}}[r]\ar^{\cdot x}[d]& \Omega e_p \ar^{\cdot (xy)^{l}}[r] \ar^{\cdot 1}[d] &   \Omega e_p\ar^{\cdot xy}[d]\ar^{\cdot (xy)^{p-l}}[r]&  \cdots   \ar^{\cdot (xy)^{l}}[r]&\Omega e_{p}\ar^{\cdot xy}[d]&
\Omega e_{l}  \ar^{\cdot x^{p-l}}[r]\ar^{\cdot y}[d]& \Omega e_p \ar^{\cdot (xy)^{l}}[r] \ar^{\cdot xy}[d] &   \Omega e_p\ar^{\cdot 1}[d]\ar^{\cdot (xy)^{p-l}}[r]&  \cdots  \ar^{\cdot (xy)^{l}}[r]&\Omega e_{p}\ar^{\cdot 1}[d]\\
\Omega e_{l+1}  \ar^{\cdot x^{p-l-1}}[r]& \Omega e_p\ar^{\cdot (xy)^{l+1}}[r] &  \Omega e_p \ar^{\cdot (xy)^{p-l-1}}[r]& \cdots \ar^{\cdot (xy)^{l+1}}[r]&\Omega e_{p}&
\Omega e_{l-1}  \ar^{\cdot x^{p-l+1}}[r]& \Omega e_p \ar^{\cdot (xy)^{l-1}}[r]&  \Omega  e_p\ar^{\cdot (xy)^{p-l+1}}[r]& \cdots  \ar^{\cdot (xy)^{l-1}}[r]&\Omega e_{p}
}
\end{equation}
for $i$ even. The action of the elements $e_{p-1}xe_p$ and $e_pye_{p-1}$ is induced by the morphism of complexes
\begin{equation}\label{eq3p}\xymatrix@C=15pt{
\Omega e_{1} \ar^{\cdot x^{p-1}}[r]& \Omega e_p \ar^{\cdot xy}[r]   &   \Omega e_p \ar^{\cdot (xy)^{p-1}}[r]& \ar^{\cdot (xy)^{p-1}}[r]\cdots&\Omega e_{p}\ar^{\cdot xy}[d]&
&  &   & &\Omega e_{p}\ar^{\cdot 1}[d]\\
& & && \Omega e_{p}&
\Omega e_{1}  \ar^{\cdot x^{p-1}}[r]& \Omega e_p \ar^{\cdot xy}[r]&  \Omega e_p \ar^{\cdot (xy)^{p-1}}[r]&\cdots \ar^{\cdot (xy)^{p-1}}[r]  & \Omega e_{p}
}
\end{equation}
for $i$ odd and
\begin{equation}\label{eq4p}\xymatrix@C=15pt{
\Omega e_{p-1}  \ar^{\cdot x}[r]& \Omega e_p \ar^{\cdot (xy)^{p-1}}[r] &   \Omega e_p\ar^{\cdot xy}[r]&  \cdots   \ar^{\cdot (xy)^{p-1}}[r]&\Omega e_{p}\ar^{\cdot xy}[d]&
&  &  & &\Omega e_{p}\ar^{\cdot 1}[d]\\
& &  &&\Omega e_{p}&
\Omega e_{p-1}  \ar^{\cdot x}[r]& \Omega e_p \ar^{\cdot (xy)^{p-1}}[r]&  \Omega  e_p\ar^{\cdot xy}[r]& \cdots  \ar^{\cdot (xy)^{p-1}}[r]&\Omega e_{p}
}
\end{equation}
for $i$ even.

Given these actions, we can now explicitly describe the quasi-isomorphism between $X^i\otimes_\Omega X^{-1}$ and $X^{i-1}$.
For $l=1,\ldots,p$, $$X^i\otimes_\Omega X^{-1}e_l = X^i\otimes_\Omega \operatorname{cone}(\Omega e_p \langle -l\rangle [-l]\overset{\cdot y^{l}}{\to}\Omega e_{p-l}) \cong\operatorname{cone}( X^ie_p\langle -l\rangle [-l]\overset{\cdot y^{l}}{\to} X^ie_{p-l})$$ is the iterated cone of the total complex of the double complex (where we omit gradings for readability)
$$\xymatrix@C=15pt{ &  &   & &\Omega e_{p
}\ar^{\cdot (1)^{p-l}}[d] &
&  &  & &\Omega e_{p}\ar^{\cdot (1)^{p-l}}[d]
\\
 \Omega e_{p-l} \ar^{\cdot x^l}[r]& \Omega e_p \ar^{\cdot (xy)^{p-l}}[r]   &   \Omega e_p \ar^{\cdot (xy)^{l}}[r]& \ar^{\cdot (xy)^{p-l}}[r]\cdots&\Omega e_{p}&
 \Omega e_{l}  \ar^{\cdot x^{p-l}}[r]& \Omega e_p \ar^{\cdot (xy)^{l}}[r] &   \Omega e_p\ar^{\cdot (xy)^{p-l}}[r]&  \cdots  \ar^{\cdot (xy)^{p-l}}[r]&\Omega e_{p}\\
}$$
for $i$ even and $i$ odd respectively. The quasi-isomorphism of the total complex to the lower one shortened by the right-most term is then obvious.

For the $p$th summand, we have natural isomorphisms 
\begin{equation*}
\begin{split}
X^i\otimes_\Omega X^{-1}e_p &= X^i\otimes_\Omega \operatorname{cone}(\Omega e_p\langle -p\rangle [-p]{\to} 0)\\& \cong\operatorname{cone}( X^ie_p\langle -p\rangle [-p]{\to} 0) \\&\cong X^ie_p\langle -p\rangle[1-p] \cong X^{i-1}e_p.
\end{split}
\end{equation*}

We now define a number of bimodule homomorphisms, which we then show provide the multiplication maps between parts of $\HTct$.
 
\begin{lem}\label{homs} We have natural bimodule homomorphisms,
$$\beta: \Omega e_p \Omega \otimes_\Omega \Omega e_p \Omega \overset{\sim}{\to} \Omega^*,$$
$$\zeta_l: \Omega e_p \Omega \otimes_\Omega \Omega^* \overset{\sim}{\to} \Omega^*, \quad
\zeta_r: \Omega^* \otimes_\Omega \Omega e_p \Omega \overset{\sim}{\to} \Omega^*,$$
$$\epsilon: \Omega^* \otimes_\Omega \Omega^* \overset{\sim}{\to} \Omega^*,$$

$$\theta_l: \Omega \otimes_\Omega \Omega^* \twoheadrightarrow \Omega e_p \Omega, \quad
\theta_r: \Omega^* \otimes_\Omega \Omega\twoheadrightarrow \Omega e_p \Omega$$
$$\iota_l: \Omega \otimes_\Omega \Omega^* \rightarrow \Omega, \quad
\iota_r: \Omega^* \otimes_\Omega \Omega \rightarrow \Omega,$$

$$\nu_l: \Theta \otimes \Theta^\sigma \rightarrow \Omega^*, \quad
\nu_r: \Theta^\sigma \otimes \Theta \rightarrow \Omega^*.$$
\end{lem}

\proof
Firstly, $\beta$ is nothing but the bimodule isomorphism constructed in \eqref{tiltingsquare}.

The dual of the short exact sequence
$$0\to\Omega e_p\Omega\to\Omega\to \Theta \to 0$$
is
isomorphic to
\begin{equation}\label{dualseq}0 \to  \Theta^\sigma \to   \Omega^* \to \Omega e_p\Omega\to 0
\end{equation}
using the bimodule isomorphisms \eqref{omegaiso} and \eqref{thetaiso}. Applying the right exact functor unctor $\Omega e_p\Omega \otimes_\Omega -$ to \eqref{dualseq}, we obtain an exact sequence
$$\Omega e_p\Omega \otimes_\Omega\Theta^\sigma \to  \Omega e_p\Omega \otimes_\Omega \Omega^* \to \Omega e_p\Omega \otimes_\Omega\Omega e_p\Omega\to 0
$$
and noting that $\Omega e_p\Omega \otimes_\Omega\Theta^\sigma =0$, the second map is an isomorphism. The map $\zeta_l$ is then the composition $$\Omega e_p\Omega \otimes_\Omega \Omega^* \overset{\sim}{\to} \Omega e_p\Omega \otimes_\Omega\Omega e_p\Omega \overset{\beta}{\to}  \Omega^*$$ of this isomorphism with $\beta$.

Similarly $\zeta_r$ is the composition 
$$ \Omega^*\otimes_\Omega \Omega e_p\Omega  \overset{\sim}{\to} \Omega e_p\Omega\otimes_\Omega \Omega e_p\Omega \overset{\beta}{\to}  \Omega^* $$
of the isomorphism obtained by applying the right exact functor $-\otimes_\Omega\Omega e_p\Omega  $ to \eqref{dualseq} (noting that again $\Theta^\sigma \otimes_\Omega\Omega e_p\Omega =0$) with the bimodule isomorphism $\beta$.

Applying $-\otimes_\Omega \Omega^*$ to \eqref{dualseq} gives an exact sequence 
$$ \Theta^\sigma \otimes_\Omega \Omega^*\to   \Omega^* \otimes_\Omega \Omega^*\to \Omega e_p\Omega\otimes_\Omega \Omega^*\to 0$$
and, noting that $\Omega^*$ is a quotient of $(\Omega e_p)^{\bigoplus p}$ and hence $\Theta^\sigma \otimes_\Omega \Omega^*=0$, the second map is again an isomorphism. The map $\epsilon$ is the composition 
$$  \Omega^* \otimes_\Omega \Omega^*\overset{\sim}{\to}  \Omega e_p\Omega\otimes_\Omega \Omega^* \overset{\zeta_l}{\to} \Omega^*$$
of this isomorphism with the isomorphism $\zeta_l$.

The morphisms $\theta_l, \theta_r$ are just given by the compositions $$ \Omega \otimes_{\Omega} \Omega^* \cong \Omega^* \twoheadrightarrow \Omega e_p \Omega\qquad \hbox{and} \qquad \Omega^* \otimes_{\Omega} \Omega\cong \Omega^*  \twoheadrightarrow \Omega e_p \Omega$$
of the quotient map $\Omega^*  \twoheadrightarrow \Omega e_p \Omega$ from \eqref{dualseq} with the canonical isomorphisms.

We define $\iota_l, \iota_r$ as the compositions 
$$
 \Omega \otimes_{\Omega} \Omega^*\overset{\theta_l}{\to} \Omega e_p\Omega \hookrightarrow \Omega
$$
respectively
$$ 
\Omega^* \otimes_{\Omega} \Omega\overset{\theta_r}{\to} \Omega e_p\Omega \hookrightarrow \Omega
$$
of $\theta_l, \theta_r$ with the natural embedding  respectively.

The morphisms $\nu_l$ and $\nu_r$ are defined as the compositions 
$$\Theta\otimes_{\Omega}\Theta^\sigma  \to \Theta^\sigma  \hookrightarrow  \Omega^*   \qquad \hbox{and} \qquad     \Theta^\sigma\otimes_{\Omega}\Theta \to\Theta^\sigma\hookrightarrow  \Omega^* $$
of the natural actions with the embedding from \eqref{dualseq} respectively.
\endproof

To describe the product on $\HTct$ using our natural bimodule homomorphisms we split the algebra into five parts:
\begin{itemize}
\item $\Omega_-$ consisting of all copies of $\Omega$ in $\HTct^-$ (with possible shifts $\Omega\langle-lp\rangle[l(1-p)]$), 
\item $\Theta_-$, 
consisting of all copies of $\Theta$ or $\Theta^{\sigma}$ in $\HTct^-$ (of the form $\Theta^{(\sigma)}\langle-lp\rangle[l(1-p)]$),  
\item $T:=\Omega e_p \Omega\langle p\rangle[p-1]$, ,  
\item $\Theta_+$ consisting of all copies of $\Theta$ or $\Theta^{\sigma}$ in $\HTct^+$ (of the form $\Theta^{(\sigma)}\langle lp\rangle[l(p-1)-1]$),
and 
\item 
$\Omega^*_+$, consisting of all copies of  $\Omega^*$ in $\HTct^+$,(with possible shifts $\Omega^*\langle 2+lp\rangle[l(p-1)]$).
\end{itemize}

To ease checking of vanishing of multiplication due to degree reasons, we now provide a table describing in which degrees each of the $\Omega$-$\Omega$-bimodule components are concentrated. Here the first element in each list is the degree of the generators, so $j$-degrees grow successively more negative, and $k$-degrees grows successively more positive.

\begin{tabular}{c|c|c}
& nonzero $j$-degrees &nonzero $k$-degrees\\
\hline
$\Omega\langle-lp\rangle[l(1-p)]$& $-lp, \cdots, -(l+2)p+2$& $lp-l, \cdots, (l+2)(p-1)$\\
&&\\
$\Theta^{(\sigma)}\langle-lp\rangle[l(1-p)]$&$-lp,\cdots, -(l+1)p+2 $&$l(p-1),\cdots, (l+1)(p-1)-1 $\\&&\\
$\Omega e_p\Omega\langle p\rangle[p-1]$&$p,\cdots, -p+2$& $1-p,\cdots, p-1 $\\&&\\
$\Theta^{(\sigma)}\langle lp\rangle[l(p-1)-1]$&$lp,\cdots, (l-1)p+2 $&$ l(1-p)+1, \cdots, (l-1)(1-p)$\\&&\\
$\Omega^*\langle 2+lp\rangle[l(p-1)]$&$(l+2)p,\cdots,  2+lp $&$(l+2)(1-p),\cdots,  l(1-p)$.
\end{tabular}

With this information, we can now prove the following proposition.

\begin{prop}\label{clubsuitmult}
The multiplication between these five parts is given by the following table:
$$
\xymatrix@C=10pt@R=3pt{
                   & \Omega_-& \Theta_-   & T  & \Theta_+ & \Omega^*_+ \\
\Omega_-           &     a   &      a     &    a               &   0               &     \iota, \theta  , a   \\
\Theta_-           &     a   &      a     &    0    &   0, a,\nu   &       0    \\
T  &    a    &     0        &    \beta     &   0   &          \zeta         \\
\Theta_+           &    0    &  0, a,\nu  &          0       &   0    &  0  \\
\Omega^*_+         &  \iota, \theta  , a     &   0    &       \zeta          & 0 &    \epsilon    \\
}$$
Here $a$ is our generic notation for an action map. For the products where we give several options, the choice depends on the component in which the product lands. In the case of products between $\Omega_-$ and $\Omega^*_+$ this is determined by
$$\xymatrix@C=10pt@R=3pt{
\textrm{Component in which the product lands:} & \Omega_- & T & \Omega^*_+ \\
\textrm{Natural map describing the product:} & \iota & \theta & a
}$$
and in the case of products between $\Theta_-$ and $\Theta_+$, it is given by
$$\xymatrix@C=10pt@R=3pt{
\textrm{Component in which the product lands:} & \HTct^- & T & \Theta_+ & \Omega^*_+ \\
\textrm{Natural map describing the product:} & 0 &  0 & a & \nu
.}$$ 
\end{prop}

\proof
The fact that the product on $\HTct^-$ is as given in the top left $2 \times 2$-corner of our table we have already established in a previous paper \cite[Theorem 32]{MT3}. 

Thanks to our simple preserving duality, we can rephrase everything in terms of right modules (obtaining a quasi-isomorphism between  $\bft^{!-1}$ and $$Y^{-1} =\operatorname{cone}\left((e_p\Omega)^p \to \sum_{l=1}^{p-1}e_l\Omega\right)$$ and between $\bft^{!i}$ and $$Y^i =\operatorname{cone}\left(\sum_{l=1}^{p-1}e_l\Omega \to (e_p\Omega)^{p-1}\to\cdots \to (e_p\Omega)^{p-1}\to (e_p\Omega)^p\right)$$ respectively (with analogous actions to those given in \eqref{eq1},\eqref{eq2},\eqref{eq1p},\eqref{eq2p},\eqref{eq3}, \eqref{eq4},\eqref{eq3p} and \eqref{eq4p}), where we obtain an obvious quasi-isomorphism $Y^{-1} \otimes_\Omega Y^i =Y^{i-1}$, implying that it suffices to check multiplications in one order.

We next consider the bottom right $3 \times 3$ corner, which provides the multiplication on $\HTct^+$. 

Note that, $\Omega e_p \Omega$ being the tilting bimodule and quasi-isomorphic to $\bft^!$, the tensor algebra $\bbT_\Omega\Omega e_p \Omega$ is necessarily a subalgebra of $\HTct^+$. Thanks to the isomorphism $\Omega e_p \Omega \otimes_\Omega \Omega e_p \Omega \overset{\sim}{\to} \Omega^*$, the multiplicative structure of this is given by $\beta,\zeta_l, \zeta_r$ and $\epsilon$, providing the nonzero entries in this square.

The product between $\Theta_+$ and $\Theta_+$ is zero by degree reasons. Indeed, the tensor product of $\Theta^{(\sigma)}\langle lp\rangle[l(p-1)-1]$ appearing in $\bbH(\bft^{!i}) $ and $\Theta^{(\sigma)}\langle l'p\rangle[l'(p-1)-1]$ appearing in $\bbH(\bft^{!i'}) $ is generated in $j$-degree $(l+l')p$ and $k$-degree $(l+l')(1-p)+2$.  The only nonzero component of  $\bbH(\bft^{!(i+i')}) $ in this $j$-degree is the top of $\Theta^{(\sigma)}\langle (l+ l')p\rangle[(l+l')(p-1)-1]$, but this has incorrect $k$-degree.

Both ${}_\Omega \Omega^*$ and $\Omega e_p \Omega$ are quotients of $\Omega e_p^{\oplus p}$ (and using the simple-preserving duality on $\Omega$, similarly $\Omega^*_\Omega $ and $\Omega e_p \Omega$ are quotients of $(e_p\Omega)^{\oplus p}$), and $\Theta\otimes_\Omega\Omega e_p = \Theta^\sigma\otimes_\Omega\Omega e_p =0$ (and similarly $e_p\Omega\otimes_\Omega\Theta = e_p\Omega\otimes_\Omega\Theta^\sigma = 0$), thus the remaining zeros in this square follow from right exactness of $\Theta^{(\sigma)} \otimes_\Omega -$ (respectively $-\otimes_\Omega\Theta^{(\sigma)}$).

It remains to confirm the bottom left $3\times 2$  (or equivalently, top right $2\times 3$) rectangle of our table.

Repeating the argument about $\Omega^*$ and $\Omega e_p \Omega$ being quotients of sums of the $p$th projective, we obtain that mutliplications between $\Theta_-$ and $\Omega^*$ respectively $\Omega e_p \Omega$ in either order are again zero.

The fact that multiplication between $\mathbb{H}(\bft^{!}) \cong \Omega e_p\Omega\langle p \rangle [p-1]$ and $\Omega_-$ is just the normal action map follows immediately from the quasi-isomorphism between $\bft^{!}$ and $\Omega e_p\Omega$.

If the product between $\Theta_+$ and $\Theta_-$ (in either order) lands in $\HTct^-$ or $\Omega e_p \Omega$, it is again zero by degree reasons.  Indeed, the tensor product of $\Theta^{(\sigma)}\langle lp\rangle[l(p-1)-1]$ appearing in $\bbH(\bft^{!i}) $ and $\Theta^{(\sigma)}\langle -l'p\rangle[l'(1-p)]$ appearing in $\bbH(\bft^{!-i'}) $ is generated in $j$-degree $(l-l')p$ and $k$-degree $(l-l')(1-p)+1$. 
Since by assumption $i'>i$, the only subspace with this nonzero $j$ degree in $\bbH(\bft^{!i-i'})$ is  the top of $\Theta^{(\sigma)}\langle -(l'-l)p\rangle[(l'-l)(1-p)]$, but this again has the wrong $k$-degree.

 For products involving $\Omega_-$ and $\Omega^*_+$, note that $\Omega^*_+$ is a component of the subalgebra $\bbT_\Omega\Omega e_p \Omega$. Multiplications being induced by the  action maps hence follows from the same claim for $\Omega e_p \Omega$.

In order to analyse the remaining multiplications, note that thanks to \cite[Theorem 32]{MT3}, which proves that $\HTct^-$ is indeed just a tensor algebra, it suffices to consider the case where one is a component of $\mathbb{H}(\bft^{!-1})$ and the other a component of 
$\mathbb{H}(\bft^{!i})$ for $i>1$, so consider multiplication $\mathbb{H}(\bft^{!i})\otimes\mathbb{H}(\bft^{!-1})\to \mathbb{H}(\bft^{!i-1})$ coming from the quasi-isomorphism $X^i\otimes_\Omega X^{-1} \to X^{i-1}$ described before Lemma \ref{homs}.

Then products between $\Theta_+$ and $\Theta_-$ being as stated follows from $\Theta_-$ appearing as a quotient of $\bigoplus_{l=1}^{p-1}\Omega e_l$ in $X^{-1}$, the explicit maps, given in \eqref{eq1},\eqref{eq2},\eqref{eq3} and \eqref{eq4}, describing the right action of $\Omega$ on $X^{i}$, and the explicit description of how elements in terms of $X_i$ correspond to elements in $\Theta_+$ following \eqref{plushomiso}.

In order to verify that the product between $\Theta_+$ and $\Omega_-$, we again look at the explicit action maps. Indeed, since $\Omega_-$ appears as a submodule of $(\Omega e_p)^{\oplus p}$ in $X^{-1}$, a lift of an element in $\Omega_-$ to $X^{-1}$ is necessarily of the form $e_l\omega e_p$ for some $\omega$. Since in the right action of $e_l\omega e_p$ on $X^{i}$, any lift of $\Theta_+$ in $X^{i}$ is annihilated, the product between $\Theta_+$ and $\Omega_-$ is zero as stated.
\endproof

\section{Explicit Hochschild cohomology of some bimodules.} \label{explicitbimod}

Here we describe the components of 
$\HH(\bfc^!, \HTct)$ as $\HH(\bfc^!)$-$\HH(\bfc^!)$-bimodules.

We fix the element $z:= \sum_{l=2}^p xye_l$ in $\Omega$.

Let us first describe the centres of our algebras $\bfc$ and $\bfc^!$.

\begin{lem}
The centre of $\bfc$ is $Z(\bfc) = F.1 \oplus \bfc^2=\sum_{l=1}^{p-1} F\cdot \xi\eta e_l$.
The centre of $\Omega$ is $Z(\Omega) = F[z]/z^p$ where $z = xy$ has $k$-degree $2$.
\end{lem}

\begin{prop} \label{bimodulecomputation}
Suppose $p>2$.
\begin{enumerate}[(i)]
\item \label{11i} $\HH(\Omega)$ is isomorphic to $Z(\bfc) \otimes Z(\Omega) \otimes \bigwedge(\kappa) / ( \bfc^2.z, \bfc^2 \kappa, z^{p-1} \kappa)$,
where  $\bfc^2$ has $jk$-degree $(2,0)$, the $z$ has $jk$-degree $(-2,2)$ and  $\kappa$ has $jk$-degree $(0,1)$.
\item \label{11ii} $\HH(\Omega,\Theta)$ is isomorphic to $\HH(\Omega)/(z^{\frac{p-1}{2}})$ as an $\HH(\Omega)$-$\HH(\Omega)$-bimodule.
\item \label{11iii} $\HH(\Omega,\Theta^\sigma)$ is isomorphic to $\HH(\Omega, \Theta)^* \langle 4-p\rangle [2-p]$
as an $\HH(\Omega)$-$\HH(\Omega)$-bimodule.
\item\label{11iv} $\HH(\Omega,\Omega^*)$ is isomorphic to $\Omega^0$, the degree $0$ part of $\Omega$.
\item\label{11v} $\HH(\Omega,\Omega e_p \Omega)$ is isomorphic to the kernel of the natural surjection
$$\HH(\Omega) \rightarrow \HH(\Omega)/(z^{\frac{p-1}{2}}).$$
\end{enumerate}
\end{prop}

\proof

\eqref{11i} By Theorem \ref{bimodulehochschild}, we need to compute the homology of
$D_{\bfc} := \bigoplus_{s,t} e_s \bfc e_t \otimes e_t \Omega e_s$ with differential
sending $\alpha \otimes a$ to
$$\alpha \xi\otimes ya  + \alpha \eta \otimes xa -(-1)^{|a|} \xi\alpha  \otimes ay -(-1)^{|a|}\eta\alpha \otimes ax.$$
The complex $D_{\bfc}$ is $\mathbb{Z}^2$-graded, where we give $e_s$ degree $(0,0)$,
we give $x$ and $y$ degree $(0,1)$, and we give $\xi$ and $\eta$ degree $(-1,0)$.  The differential therefore has degree $(-1,1)$. We remark that this is \emph{ not} our usual $(j,k)$-grading and we still denote by $|\cdot |$ the $k$-degree of an element as before.
We have a basis for $e_s \bfc e_t \otimes e_t \Omega e_s$ given by those monomials
$e_s \xi^{m_\xi} \eta^{m_\eta} e_t \otimes e_t x^{m_x} y^{m_y}e_s$ which are not zero in this space. We set $$a_{s,l}=e_s\xi e_{s+1}\otimes e_{s+1}yz^le_s, \quad b_{s,l}=e_s\eta e_{s-1}\otimes e_{s-1}xz^le_s, \quad w_{s,l}= e_s \xi\eta e_s \otimes e_s z^{l+1}e_s$$ 
and note that $a_{s,l}\neq 0$ if and only if $l+1 \leq s \leq p-1$, $b_{s,l}\neq 0$ if and only if $l+2\leq s\leq p$ and $w_{s,l}\neq 0$ if and only $l+2 \leq s\leq p-1$. Moreover, $a_{s,l}$ and $b_{s,l}$and $w_{s,l}$ vanish for all $s$ if $ l \geq p-1$. The nonzero graded subspaces of $D_{\bfc}$ are $D_{\bfc}^{-2,0}$ (which is  just $\bfc^2 \otimes 1_\Omega$ and isomorphic to $\bfc^2$), $D_{\bfc}^{0,2l}$ for $0\leq l \leq p-1$, $D_{\bfc}^{-1,2l+1}$ and $D_{\bfc}^{-2,2l+2}$ for $0 \leq l \leq p-2$. The first is just $\bfc^2 \otimes 1_\Omega$ and isomorphic to $\bfc^2$. For fixed $l$, $D_{\bfc}^{0,2l}$, $D_{\bfc}^{-1,2l+1}$ and $D_{\bfc}^{-2,2l+2}$ have  bases given by $\{e_s\otimes e_s z^l e_s | s= l+1, \dots p\}$, $\{a_{s,l} |  s= l+1, \dots p-1\}\cup\{b_{s,l} |  s= l+2, \dots p\}$ and $\{w_{s,l}| s=l+2,\dots,p-1\}$ respectively.
Our complex $D_{\bfc}$ is then a sum of the complex
$$0 \rightarrow \bfc^2 \rightarrow 0$$
and the sum over $l$ of complexes, for $0 \leq l \leq p-2$,
\begin{equation*}
\begin{split}
(0 \rightarrow &D_{\bfc}^{(0,2l)} \rightarrow D_{\bfc}^{(-1,2l+1)} \rightarrow D_{\bfc}^{(-2,2l+2)} \rightarrow 0)\\ &\cong(0 \rightarrow F^{p-l} \rightarrow F^{2p-2-2l} \rightarrow F^{p-2-l} \rightarrow 0)
\end{split}
\end{equation*}
(where we interpret spaces as zero if they have zero or negative dimensions, which happens for $D_{\bfc}^{(-2,2l+2)} $ for $l\geq p-2$ and for $D_{\bfc}^{(-1,2l+1)}$ for $l=p-1$)
and the differential  acts  on the $l$-component by
\begin{equation*}\begin{split}
e_s\otimes e_s z^l e_s &\mapsto a_{s,l}+b_{s,l}-a_{s-1,l}-b_{s+1,l},\\
 a_{s,l} &\mapsto w_{s,l}-w_{s+1,l}\\
b_{s,l} &\mapsto -w_{s,l}+ w_{s-1,l}
\end{split}\end{equation*}
from where we see that in the sequence $D_{\bfc}^{0,2l}\to D_{\bfc}^{-1,2l+1} \to D_{\bfc}^{-2,2l+2}$ the last map is surjective, the first has one-dimensional kernel spanned by $\sum_{s=l+1}^{p} e_s\otimes e_s z^l e_s = 1\otimes z^l$ (which lies in the centre of $\Omega$), and one-dimensional homology in the middle spanned by $\kappa z^l$ where $\kappa:= \sum_{s=1}^{p-1} a_{s,0}$.  The homology $\bbH(D_{\bfc})$ is therefore
$$c^2 \oplus \bigoplus_{l=0}^{p-2} F.\kappa z^l \oplus \bigoplus_{l = 0}^{p-1} F.z^l$$ and the multiplication is obvious from this explicit description and Proposition \ref{cup}.
 In our gradings, the $j$-grading sees $\eta, \xi, x,y$ in degrees $1,1,-1,-1$ respectively, and the $k$ grading has  $\eta, \xi, x,y$ in degrees $0,0,1,1$, so the factor $\bfc^2$ has $(j,k)$-degree $(2,0)$, the element $z$ has $(j,k)$-degree $(-2,2)$ and the element $\kappa$ has $(j,k)$-degree $(0,1)$. This completes the proof of \eqref{11i}.

\eqref{11ii} By Theorem \ref{bimodulehochschild}, we need to compute the homology of $D_{\bfc,\Theta}:=\bigoplus_{s,t} e_s \bfc e_t \otimes e_t \Theta e_s$ with differential
$$\alpha \otimes m \mapsto  \alpha \xi\otimes ym  + \alpha \eta \otimes xm -(-1)^{|m|} \xi\alpha  \otimes my -(-1)^{|m|}\eta\alpha \otimes mx.$$
Using the same grading and notation as in \eqref{11i}, the only nonzero graded components are $D_{\bfc,\Theta}^{-2,0}$ (which, as before, is just $\bfc^2 \otimes 1_\Omega\cong \bfc^2$ and contributes to homology), $D_{\bfc,\Theta}^{0,2l}$, $D_{\bfc,\Theta}^{-1,2l+1}$ for $0 \leq l \leq \frac{p-3}{2}$ (recall that $p$ is odd) and $D_{\bfc,\Theta}^{-2,2l+2}$ for $0 \leq l \leq \frac{p-5}{2}$. When nonzero, the spaces $D_{\bfc,\Theta}^{0,2l}$, $D_{\bfc,\Theta}^{-1,2l+1}$ and $D_{\bfc,\Theta}^{-2,2l+2}$  have bases given by $\{e_s\otimes e_s z^l e_s | s= l+1, \dots p-l-1\}$, $\{a_{s,l}|  s= l+1, \dots p-l-2\}\cup\{b_{s,l} |  s= l+2, \dots p-l-1\}$ and $\{w_{s,l}| s=l+2,\dots,p-l-2\}$ respectively.
Our complex $D_{\bfc,\Theta}$ is then a sum of the complex
$$0 \rightarrow \bfc^2 \rightarrow 0$$
and the sum over $l$ of complexes  for $0 \leq l \leq \frac{p-3}{2}$,
\begin{equation*}
\begin{split}
(0 \rightarrow &D_{\bfc,\Theta}^{(0,2l)} \rightarrow D_{\bfc,\Theta}^{(-1,2l+1)} \rightarrow D_{\bfc,\Theta}^{(-2,2l+2)} \rightarrow 0)\\ &\cong(0 \rightarrow F^{p-2l-1} \rightarrow F^{2p-4l-4} \rightarrow F^{p-2l-3} \rightarrow 0)
\end{split}
\end{equation*}
and the differential acts as before on the basis elements. Again the last map is surjective, the first has kernel $\sum_{s=l+1}^{p-l-1} e_s\otimes e_s z^l e_s = 1\otimes z^l$, and homology in the middle is spanned by $\kappa z^l= \sum_{s=l+1}^{p-l-2} a_{s,l}$.
The homology $\bbH(D_{\bfc,\Theta})$ is therefore
$$\bfc^2 \oplus \bigoplus_{l=0}^{\frac{p-3}{2}} F.z^l \kappa \oplus \bigoplus_{l = 0}^{\frac{p-3}{2}} F.z^l. $$
By Proposition \ref{cup}, the $\HH(\Omega)$-$\HH(\Omega)$-bimodule structure induced by multiplication in $\bfc^{\op}$ and the $\Omega$-$\Omega$-bimodule structure on $\Theta$. Using the explicit description of basis elements in terms of tensor products of elements in $\bfc$ and elements in $\Theta$, \eqref{11ii} follows.

\eqref{11iii} Again by Theorem \ref{bimodulehochschild}, in order to compute $\HH(\Omega, \Theta^\sigma)$ 
we need to compute the homology of $$D_{\bfc,\Theta^\sigma}:= \bigoplus_{s,t} e_s \bfc e_t \otimes e_t \Theta^\sigma e_s,$$ 
with differential
\begin{equation*}\begin{split}
\alpha \otimes m \mapsto & \alpha \xi\otimes ym  + \alpha \eta \otimes xm -(-1)^{|m|} \xi\alpha  \otimes m\cdot y -(-1)^{|m|}\eta\alpha \otimes m\cdot x \\&=  \alpha \xi\otimes ym  + \alpha \eta \otimes xm -(-1)^{|m|} \xi\alpha  \otimes mx -(-1)^{|m|}\eta\alpha \otimes my 
\end{split}\end{equation*}
where we denote by $m\cdot x$ the action of $x\in \Omega$ on the element $m \in \Theta^\sigma$ and by $mx$ the usual (untwisted) action of $\Omega$ on $\Theta$.
As a vector space, this is isomorphic to $\bigoplus_{s,t} e_s \bfc e_t \otimes e_t \Theta e_{p-s}$. This has nonzero components
$D_{\bfc,\Theta^\sigma}^{(0, p-2-2l)}$ for $l=0, \dots, \frac{p-3}{2}$, as well as $D_{\bfc,\Theta^\sigma}^{(-1, p-1-2l)}$ and $D_{\bfc,\Theta^\sigma}^{(-2, p-2l)}$ for $l=1, \dots, \frac{p-1}{2}$, with bases given by
$$\{e_s \otimes e_sx^{p-s-l-1}y^{s-l-1}e_{p-s}  | s=l+1,\dots,p-l-1\}$$
$$\{e_s\xi e_{s+1} \otimes e_{s+1}x^{p-s-l-1}y^{s-l}e_{p-s}, e_{s+1}\eta e_{s} \otimes e_{s}x^{p-s-l-1}y^{s-l}e_{p-s-1}  | s=l,\dots,p-l-1\}$$
and
$$\{e_s \xi \eta e_s\otimes e_sx^{p-s-l}y^{s-l}e_{p-s}  | s=l,\dots,p-l\}$$
respectively.
As the differential has degree $(-1,1)$, for $l=0$ we obtain homology spanned by $\{e_s \otimes e_sx^{p-s-1}y^{s-1}e_{p-s}  | s=1,\dots,p-1\}$ in degree $(0,p-2)$. This is equal to $1 \otimes (\Theta^\sigma)^{p-2}$. The rest of the complex is a sum over $l$ for $l=1, \dots, \frac{p-1}{2} $ of
\begin{equation*}
\begin{split}
(0 \rightarrow &D_{\bfc,\Theta^\sigma}^{(0,p-2-2l)} \rightarrow D_{\bfc,\Theta^\sigma}^{(-1,p-1-2l)} \rightarrow D_{\bfc,\Theta^\sigma}^{(-2,p-2l)} \rightarrow 0)\\ &\cong(0 \rightarrow F^{p-2l-1} \rightarrow F^{2p-4l} \rightarrow F^{p-2l+1} \rightarrow 0).
\end{split}
\end{equation*}
Setting $$f_{s,l} = e_s\xi e_{s+1} \otimes e_{s+1}x^{p-s-l-1}y^{s-l}e_{p-s}, $$ $$ g_{s,l} = e_{s+1}\eta e_{s} \otimes e_{s}x^{p-s-l-1}y^{s-l}e_{p-s-1}$$ $$ v_{s,l} = e_s \xi \eta e_s\otimes e_sx^{p-s-l}y^{s-l}e_{p-s}$$ respectively, the differential acts as
\begin{equation*}
\begin{split}e_s \otimes e_sx^{p-s-l-1}y^{s-l-1}e_{p-s} &\mapsto f_{s,l} + f_{s-1,l} + g_{s,l} + g_{s-1,l}\\
f_{s,l} & \mapsto v_{s,l} + v_{s+1,l}\\ g_{s,l} &\mapsto - v_{s,l} - v_{s+1,l}.
\end{split}
\end{equation*}
It is easy to see that the first map is injective. However, the image of the last map is spanned by $v_{s,l}+v_{s+1,l}$ for $s = l, \dots, p-l-1$ and is hence only $p-2l$-dimensional, leaving one-dimensional homology in both the middle (spanned by $\mu_l = (f_{\frac{p-1}{2},l} + g_{\frac{p-1}{2},l})$ say) and the end (spanned by $\nu_l = (v_{\frac{p-1}{2},l} - v_{\frac{p+1}{2},l})$, say).
In order to describe the structure as $\HH(\Omega)$-$\HH(\Omega)$-bimodule, we need to determine the action of the generators of $\HH(\Omega)$ on this, and in light of Proposition \ref{cup} this is induced by multiplication in $\bfc^{\op}$ and the action of $\Omega$ on either side of $\Theta^\sigma$, or, in other words, the natural action of $D_\bfc$ on $D_{\bfc,\Theta^\sigma}$. It is clear that both $\nu_l$ and $\mu_l$ are annihilated by $\bfc^2$. Direct computation shows that $\kappa.\mu_l = \mu_l.\kappa= \frac{1}{2}\nu_l$, $z.\mu_l = \mu_l.z=\mu_{l-1}$ and $z.\nu_l=\nu_l.z = \nu_{l-1}$.
By graded dimensions, the only other non-zero product could be $\bfc^2 . D_{\bfc,\Theta^\sigma}^{(0, p-2)}$, which lies in degree $(-2,p-2)$, where $\nu_1$ also lives. Direct computation shows that with our choice of representatives of homology, we obtain 
\begin{equation*}
\begin{split}
(e_s\xi\eta e_s \otimes e_s)(e_s \otimes e_s x^{p-s-1}y^{s-1}e_{p-s})& 
= (e_s \otimes e_s x^{p-s-1}y^{s-1}e_{p-s})(e_s\xi\eta e_s \otimes e_s)\\&=\frac{1}{2}(-1)^{\frac{p-1}{2} -s} \nu_1
\end{split}
\end{equation*}
 and all other product with non-matching idempotents are obviously zero.
The $(j,k)$-degrees of the basis elements are $(-p+2, p-2)$ for $e_s \otimes e_sx^{p-s-1}y^{s-1}e_{p-s} $ for $s=1,\dots,p-1$, then $(-p+2+2l,p-2l-1)$ for $\mu_l$ and $(-p+2+2l,p-2l)$ for $\nu_l$.

This completes our combinatorial description of $\HH(\Omega, \Theta^\sigma)$.
To define an isomorphism between $\HH(\Omega,\Theta^\sigma)$ and $\HH(\Omega, \Theta)^*$ we now define a bilinear form
$$\lvert -,- \rvert: \HH(\Omega,\Theta^\sigma) \otimes \HH(\Omega, \Theta) \rightarrow F$$
such that 
$$\lvert h,h'h'' \rvert =\lvert hh',h'' \rvert, \qquad \lvert h,h''h' \rvert = (-1)^{|h'|_k(|h|_k+|h''|_k)} \lvert h'h,h'' \rvert,$$
for $h \in \HH(\Omega, \Theta^\sigma)$, $h' \in \HH(\Omega)$, $h'' \in \HH(\Omega, \Theta)$.
Indeed the form $\lvert -,- \rvert$ which pairs $2(-1)^{\frac{p-1}{2} -s}(e_s \otimes e_sx^{p-s-1}y^{s-1}e_{p-s}) \in  \HH(\Omega,\Theta^\sigma)$ (of $(j,k)$-degree $(2-p,p-2)$) 
with $e_s\xi\eta \otimes 1 \in  \HH(\Omega,\Theta)$ (which has $(j,k)$-degree $(2,0)$), 
which pairs $z^l$ (of $(j,k)$-degree $(-2l,2l)$) with $\nu_{l+1}$ (of $(j,k)$-degree $(2l-p+4,p-2l-2)$), and which pairs $z^l\kappa$ (of $(j,k)$-degree $(-2l,2l+1)$) with $\mu_{l+1} $(of $(j,k)$-degree $(2l-p+4,p-2l-3)$) has the required property; 
in fact all signs $(-1)^{|h'|_k(|h|_k+|h''|_k)}$ are $+1$ when $\lvert h'h,h'' \rvert$ is nonzero for elements $h$,$h'$,$h''$ 
of our canonical bases
since the super-commutation relations defining $\HH(\Omega)$ are all commutation relations, with $z$ lying in degree $2$. It follows that there is an isomorphism $$\HH(\Omega,\Theta^\sigma) \cong\HH(\Omega, \Theta)^* \langle 4-p\rangle [2-p]$$ as claimed.

\eqref{11iv} Similarly to the previous ones, we apply Theorem \ref{bimodulehochschild} and see that we need to compute the homology of the complex $D_{\bfc, \Omega^*}:= \bigoplus_{s,t} e_s \bfc e_t \otimes e_t \Omega^* e_s$ with differential
sending $\alpha \otimes \varphi$ to
$$\alpha \xi\otimes y\varphi  + \alpha \eta \otimes x\varphi -(-1)^{|\varphi|} \xi\alpha  \otimes \varphi y -(-1)^{|a|}\eta\alpha \otimes \varphi x.$$

The computation is similar to the one in \eqref{11i}. We set  $$\tilde{a}_{s,l-1}=e_{s-1} \xi e_s\otimes e_s(z^{l-1}x)^*e_{s-1} \qquad \quad \tilde{b}_{s,l-1}=e_{s+1}\eta e_s\otimes e_{s}(z^{l-1}y)^*e_{s+1}.$$ The nonzero graded components of $D_{\bfc, \Omega^*}$ are $D_{\bfc, \Omega^*}^{(0,0)}$, having basis given by $\{e_s\otimes e_s^* \mid s=1,\dots p\}$, as well as $D_{\bfc, \Omega^*}^{(0,-2l)}, D_{\bfc, \Omega^*}^{(-1,-2l+1)}$ and $D_{\bfc, \Omega^*}^{(-2,-2l+2)} $ for $1\leq l \leq p-1$ with respective bases 
given by $$\{e_s\otimes e_s(z^l)^*e_s \mid s=l+1, \dots, p\}$$
$$\{e_{s-1} \xi e_s\otimes e_s(z^{l-1}x)^*e_{s-1}, e_s\eta e_{s-1}\otimes e_{s-1}(z^{l-1}y)^*e_s \mid s=l+1, \dots, p\}$$
$$\{e_s\xi\eta e_s\otimes e_s(z^{l-1})^*e_s \mid s=l, \dots, p-1\}.$$
Our complex is isomorphic to the direct sum of $p$ complexes
$$0 \rightarrow D_{\bfc, \Omega^*}^{(0,-2l)} \rightarrow D_{\bfc, \Omega^*}^{(-1,-2l+1)} \rightarrow D_{\bfc, \Omega^*}^{(-2,-2l+2)} \rightarrow 0
$$
for $l=1,..,p-1$ and $0 \rightarrow D_{\bfc, \Omega^*}^{(0,0)} \to 0$.
The last summand provides the homology claimed in this case, so we need to show that the first $p-1$ summands are exact.
Indeed, the dimensions of $D_{\bfc, \Omega^*}^{(0,-2l)}, D_{\bfc, \Omega^*}^{(-1,-2l+1)}$ and $D_{\bfc, \Omega^*}^{(-2,-2l+2)} $ are $p-l$, $2(p-l)$ and $p-l$ respectively, so it suffices to show that the differential is injective on the first and surjective on the last component.
Since $$\tilde{a}_{s,l-1}=e_{s-1} \xi e_s\otimes e_s(z^{l-1}x)^*e_{s-1} \qquad \quad \tilde{b}_{s,l-1}=e_{s+1}\eta e_s\otimes e_{s}(z^{l-1}y)^*e_{s+1},$$
the differential acts as
$$ e_s\otimes e_s(z^l)^*e_s \mapsto a_{s+1,l-1}+b_{s-1,l-1}- a_{s,l-1}-b_{s,l-1},$$
where summands are considered as zero if $s$ falls outside of the range $1,\dots, p$,
from which we see injectivity of the first differential.
The basis element $\tilde{a}_{s,l-1}$ in $D_{\bfc, \Omega^*}^{(-1,-2l+1)}$ gets sent to
$e_{s-1}\xi\eta e_{s-1}\otimes e_{s-1}(z^{l-1})^*e_{s-1} -e_s\xi\eta e_s\otimes e_s(z^{l-1})^*e_s$
where again summands are considered as zero if $s$ falls outside of the range $1,\dots, p$, from which we see surjectivity of the the second differential, completing the proof of \eqref{11iv}.

\eqref{11v} We have an exact sequence of $\Omega$-$\Omega$-bimodules,
$$0 \rightarrow \Omega e_p \Omega \rightarrow \Omega \rightarrow \Theta \rightarrow 0.$$
Applying $\RHom_{\Omega\otimes \Omega^{\op}}(\Omega,-)$ gives us an exact triangle
$$\RHom_{\Omega \otimes \Omega^{\op}}(\Omega, \Omega e_p \Omega) \rightarrow \RHom_{\Omega \otimes \Omega^{\op}}(\Omega, \Omega)
\rightarrow \RHom_{\Omega \otimes \Omega^{\op}}(\Omega,\Theta) \rightsquigarrow$$
in the derived category of $F$-$F$-bimodules, which corresponds to an exact triangle
$$\HH(\Omega, \Omega e_p \Omega) \rightarrow \HH(\Omega, \Omega)
\rightarrow \HH(\Omega,\Theta) \rightsquigarrow$$
We know $\HH(\Omega, \Omega)$ and $\HH(\Omega,\Theta)$, and from our calculations
the map between them is visibly the canonical surjection. This completes the proof of (v).
\endproof

We give some pictures visualising the structure of the bimodules in case $p=5$ (the numbers down the left hand side denote the $k$-grading and along the top the j-grading 
Here is $\HH(\Omega)$:
$$\xymatrix@C=10pt@R=3pt{ & 2&1& 0& -1& -2& -3&-4&-5&-6&-7&-8\\
0 \hspace{1cm}&     & &\ar@{-}[dll] 1 \ar@{-}[dd] \ar@{-}[dddrr]  &&&&    &          & & & \\
0 \hspace{1cm}& F^{\oplus p-1} &&         &                       & & & \\
1 \hspace{1cm}&     & &\kappa \ar@{-}[ddrr] &           & & & &&&&\\
2 \hspace{1cm}&     &         &    &                   & z  \ar@{-}[d] \ar@{-}[ddrr]      & & & \\
3 \hspace{1cm}&     &         & & &\kappa z \ar@{-}[ddrr] & & & &&&\\
4 \hspace{1cm}&     &         &                       &&& & z^2 \ar@{-}[d] \ar@{-}[ddrr] & & &&\\
5 \hspace{1cm}&     &         &       &&  & & \kappa z^2 \ar@{-}[ddrr]  &&&& \\
6 \hspace{1cm}&     &         &        &&  &            && & z^3 \ar@{-}[d] \ar@{-}[ddrr]& &\\
7 \hspace{1cm}&     &         &         & & &&           && z^3 \kappa  &&\\
8 \hspace{1cm} &     &         &          &            &            &&&&& & z^4
}$$
Here is $\HH(\Omega, \Theta)$:
$$\xymatrix@C=10pt@R=3pt{&2&1&0&-1&-2&-3&&\\
0 \hspace{1cm}&     && \ar@{-}[dll] 1 \ar@{-}[dd] \ar@{-}[dddrr]      &          & & & \\
0 \hspace{1cm}& F^{\oplus p-1} &&         &                       & & & \\
1 \hspace{1cm}&     && \kappa \ar@{-}[ddrr] &           & & & \\
2 \hspace{1cm}&     &&         &                       & z  \ar@{-}[d]      & & & \\
3 \hspace{1cm}&     &&         && \kappa z   & & & \\
}$$
Here is $\HH(\Omega, \Theta^\sigma)$:
$$\xymatrix@C=10pt@R=3pt{&1&0&-1&-2&-3&&\\
0 \hspace{1cm}&                            \mu_2 \ar@{-}[d] \ar@{-}[ddrr] &     &          & & & \\
1 \hspace{1cm}& \nu_2 \ar@{-}[dddrr]                            &&           & & & \\
2 \hspace{1cm}&                           &                        &         \mu_1  \ar@{-}[dd]      &&& &  \\
3 \hspace{1cm}&    &         &                       &                  & \ar@{-}[dll] F^{\oplus p-1} &  & \\
3 \hspace{1cm}&                           &         &  \nu_1 & & & & \\
}$$
Here is $\HH(\Omega, \Omega^*)$:
$$\xymatrix@C=10pt@R=3pt{
0 \hspace{1cm}&     &         &                       & & F^{\oplus p} & & 
}$$
Here is $\HH(\Omega, \Omega e_p\Omega)$:
$$\xymatrix@C=10pt@R=3pt{ && & & &-4&-5&-6&-7&-8\\
4 \hspace{1cm}&     &         &                       & & z^2 \ar@{-}[d] \ar@{-}[ddrr] & && &\\
5 \hspace{1cm}&     &         &          && \kappa z^2 \ar@{-}[ddrr] & & &&\\
6 \hspace{1cm}&     &         &          &            && & z^3 \ar@{-}[d] \ar@{-}[ddrr]& &\\
7 \hspace{1cm}&     &         &          &           && & z^3 \kappa & &\\
8 \hspace{1cm} &     &         &          &            &           && & & z^4
}$$

\begin{rk} The bimodule isomorphism 
$$\HH(\Omega,\Theta^\sigma) \cong \HH(\Omega, \Theta)^* \langle 4-p\rangle [2-p]$$ 
of Proposition \ref{bimodulecomputation}\eqref{11iii} is striking, since we also have $\Theta^\sigma \cong \Theta^*$ as bimodules. 
This duality between Hochschild cohomologies does not follow from basic general principles and therefore 
deserves further comment. 
We give a more conceptual explanation of its origin here.
Thanks to \eqref{omegaiso} and \eqref{thetaiso}, the dual of the short exact sequence
$$0 \rightarrow \Omega e_p \Omega \rightarrow \Omega \rightarrow \Theta \rightarrow 0$$
is isomorphic to
$$0 \leftarrow \Omega e_p \Omega \langle 2p-2 \rangle [2p-2] 
\leftarrow \Omega^* \leftarrow \Theta^\sigma \langle p-2 \rangle [p-2] \leftarrow 0.$$
Applying derived $\Hom_{\Omega \otimes \Omega^{\op}}(\Omega,-)$ gives us an exact triangle
$$\HH(\Omega, \Theta^\sigma) \langle p-2 \rangle [p-2]
\rightarrow \HH(\Omega, \Omega^*) \rightarrow \HH(\Omega, \Omega e_p \Omega) \langle 2p-2 \rangle [2p-2] \rightsquigarrow$$
We know from Proposition \ref{bimodulecomputation}\eqref{11v} that $\HH(\Omega, \Omega e_p \Omega)$ is the kernel of $\HH(\Omega, \Omega \rightarrow \Theta)$, an extension of 
$F[\kappa,z]/(\kappa^2, z^{\frac{p-1}{2}}) \langle 1-p \rangle [1-p]$ by $F \langle 2-2p \rangle [2-2p]$ and  
we know that $\HH(\Omega, \Omega^*)$ is isomorphic to $F^{\oplus p} \langle 0 \rangle [0]$. 
Two copies of $F$ cancel in the derived category in our triangle via the map $\HH(\gamma)$ where $\gamma$ is the natural surjection $\Omega^* \twoheadrightarrow   \Omega e_p \Omega$ from \eqref{dualseq}
(see proof of Lemma \ref{crazysymbols}, the product $\lozenge_l$), 
leaving us with an exact triangle
$$\HH(\Omega, \Theta^\sigma) \langle p-2 \rangle [p-2] \rightarrow F^{\oplus p-1} \rightarrow F[\kappa,z]/(\kappa^2, z^{\frac{p-1}{2}}) \langle 1-p \rangle [1-p] \langle 2p-2 \rangle [2p-2] \rightsquigarrow,$$
that is 
$$\HH(\Omega, \Theta^\sigma) \langle p-2 \rangle [p-2] \rightarrow F^{\oplus p-1} \rightarrow F[\kappa,z]/(\kappa^2, z^{\frac{p-1}{2}}) \langle p-1 \rangle [p-1] \rightsquigarrow,$$
which we can shift to a triangle
$$F[\kappa,z]/(\kappa^2, z^{\frac{p-1}{2}}) \langle 1 \rangle [0] \rightarrow 
\HH(\Omega, \Theta^\sigma) \rightarrow F^{\oplus p-1} \langle 2-p \rangle [2-p] \rightsquigarrow$$
or
$$F[\kappa,z]/(\kappa^2, z^{\frac{p-1}{2}}) \langle p-3 \rangle [p-2] \rightarrow 
\HH(\Omega, \Theta^\sigma) \langle p-4 \rangle [p-2]\rightarrow F^{\oplus p-1} \langle -2 \rangle [0] \rightsquigarrow$$

This is dual to the exact triangle 
$$F^{\oplus p-1} \langle 2 \rangle [0] \rightarrow \HH(\Omega, \Theta) 
\rightarrow F[\kappa,z]/(\kappa^2, z^{\frac{p-1}{2}}) \rightsquigarrow.$$
Here we use the self-injectivity of $F[\kappa,z]/(\kappa^2, z^{\frac{p-1}{2}})$,
which is given by an isomorphism 
$$F[\kappa,z]/(\kappa^2, z^{\frac{p-1}{2}}) \cong F[\kappa,z]/(\kappa^2, z^{\frac{p-1}{2}})^* \langle 3-p \rangle [2-p]$$
of $F[\kappa,z]/(\kappa^2, z^{\frac{p-1}{2}})$-$F[\kappa,z]/(\kappa^2, z^{\frac{p-1}{2}})$-bimodules.
We thus have
$$\HH(\Omega, \Theta)^* \cong \HH(\Omega, \Theta^\sigma) \langle p-4 \rangle [p-2]$$
as $jk$-graded $\HH(\Omega)$-$\HH(\Omega)$-bimodules.
\end{rk}

\begin{rk}\label{bases}
The spaces computed in Proposition \ref{bimodulecomputation} come with natural bases. 
Denote $\chi := \HH(\Omega)$ and $\overline{\chi} := \chi/z^{\frac{p-1}{2}}$,
and let $\underline{\chi}$ denote the kernel of the natural surjection $\chi \rightarrow \overline{\chi}$, so we have isomorphisms
$\HH(\Omega, \Theta) \cong \overline{\chi}$, $\overline{\chi}^\circledast:=\HH(\Omega, \Theta^\sigma) \cong \overline{\chi}^*\langle 4-p \rangle [2-p]$ and $\HH(\Omega,\Omega e_p \Omega) \cong \underline{\chi}$.
We have bases for these bimodules, indexed by pairs $(d,e)$ where $d$ denotes a $jk$-degree and $e$ an 
idempotent such that $em_{d,e} = m_{d,e}$ (as an example, $m_{-2l,2l,1}$ corresponds to $z^l=1\cdot z^l$,  $m_{2,0,e_s}$ corresponds to $e_s\xi\eta e_s\otimes e_s$, etc):
\begin{equation*}\begin{split}
\mathcal{B}_\chi =& \{ m_{-2l,2l,1} | 0 \leq l \leq p-1 \} \cup \{ m_{-2l,2l+1,1} | 0 \leq l \leq p-2 \}\\& \cup 
\{ m_{2,0,e_s} | 1 \leq s \leq p-1 \}; \\
\mathcal{B}_{\overline{\chi}} = & \{ m_{-2l,2l,1} | 0 \leq l \leq \frac{p-3}{2} \} \cup \{ m_{-2l,2l+1,1} | 0 \leq l \leq \frac{p-3}{2} \}
\\&\cup \{ m_{2,0,e_s} | 1 \leq s \leq p-1 \}; \\
\mathcal{B}_{\overline{\chi}^*} = & \{ m_{2l,-2l,1} | 0 \leq l \leq \frac{p-3}{2} \} \cup \{ m_{2l,-2l-1,1} | 0 \leq l \leq \frac{p-3}{2} \}\\&
\cup \{ m_{-2,0,e_s} | 1 \leq s \leq p-1 \}; \\
\mathcal{B}_{\underline{\chi}} = & \mathcal{B}_\chi \backslash \mathcal{B}_{\overline{\chi}};\\
\mathcal{B}_{\Omega^0} = & \{ m_{0,0,e_s} | 1 \leq s \leq p \}. 
\end{split} \end{equation*}
More precisely we have 
\begin{equation*}\begin{split}
\mathcal{B}_\chi = & \{ 1, z^l | 0 \leq l \leq p-1 \} \cup \{ \kappa z^l | 1 \leq l \leq p-2 \} \cup \{ e_s \xi\eta \otimes 1 | 1 \leq s \leq p-1 \}; \\
\mathcal{B}_{\overline{\chi}^\circledast} = & \{ \nu_{l+1} | 0 \leq l \leq \frac{p-3}{2}\} \cup \{ \mu_{l+1} | 0 \leq l \leq \frac{p-3}{2}\}
\\&\cup \{ e_s \otimes e_s x^{p-s-1}y^{s-1}e_{p-s}| 1 \leq s \leq p-1 \} = \{m_{j+4-p,k+p-2,e}|m_{j,k,e} \in \mathcal{B}_{\overline{\chi}^*}\}
\end{split} \end{equation*}
and we identify $\mathcal{B}_{\overline{\chi}}$ and $\mathcal{B}_{\underline{\chi}}$ with subsets of $\mathcal{B}_{\chi}$ in the natural way. 
The basis $\mathcal{B}_{\Omega^0}$ is merely the set of idempotents $e_s$ for $1 \leq s \leq p$.
\end{rk}

\section{The algebra $\HHTct = \ffHH(\HTct)$.} \label{spadesuitsection}

Cute as $\HTct$ is, to compute the Hochschild cohomology of blocks of polynomial representations of $\GL_2$ we must diminish it,
by taking Hochschild cohomology with respect to $\Omega$. The resulting algebra we call $\HHTct$.
In the remaining parts of the paper we assume $p>2$.

\subsection{Description via bimodules.} Recall the notations from Remark \ref{bases}. By taking componentwise Hochschild cohomology we see that the structure of $\HHTct$ as an ungraded $\chi$-$\chi$-bimodule is given by
$$
\xymatrix@C=6pt@R=3pt{
&&&&&& & ... & & \\
&&&&&&             & \overline{\chi}^* & \overline{\chi}          & \Omega^0 \\
&&&&&&   & \overline{\chi} & \Omega^0 & \\
&&&&&&   & \Omega^0  & & \\
&&&&&& \underline{\chi} & & \\
&&&&& \chi &&&        \\
&&&& \chi   & \overline{\chi}^*        \\
&&& \chi  & \overline{\chi}^* & \overline{\chi}       \\
&& \chi & \overline{\chi}^* & \overline{\chi} & \overline{\chi}^*        \\
& \chi  & \overline{\chi}^* & \overline{\chi} & \overline{\chi}^* & \overline{\chi}        \\
& & & ... & &
}$$

From the structure of $\HTct$ as bigraded $\Omega$-$\Omega$-bimodule, we infer the structure of $\HHTct^- = \ffHH(\HTct^-)$ as a $k$-graded $\chi$-$\chi$-bimodule
$$\xymatrix@C=6pt@R=3pt{
&&&&& \chi         \\
&&&& \chi[1-p]   & \overline{\chi}^*[2-p]        \\
&&& \chi[2-2p]  & \overline{\chi}^*[3-2p] & \overline{\chi}       \\
&& \chi[3-3p]  & \overline{\chi}^*[4-3p] & \overline{\chi}[1-p] & \overline{\chi}^*[2-p]        \\
& \chi[4-4p]  & \overline{\chi}^*[5-4p] & \overline{\chi}[2-2p] & \overline{\chi}^*[3-2p] & \overline{\chi}        \\
& & & ... & &;
}$$
the structure of $\HHTct^-$ as a $j$-graded $\chi$-$\chi$-bimodule
$$\xymatrix@C=6pt@R=3pt{
&&&&& \chi         \\
&&&& \chi \langle -p \rangle  & \overline{\chi}^*  \langle 4-p \rangle      \\
&&& \chi \langle -2p \rangle & \overline{\chi}^* \langle 4-2p \rangle  & \overline{\chi}      \\
&& \chi \langle -3p \rangle & \overline{\chi}^* \langle 4-3p \rangle  & \overline{\chi} \langle -p \rangle & \overline{\chi}^*  \langle 4-p \rangle      \\
& \chi \langle -4p \rangle & \overline{\chi}^* \langle 4-4p \rangle  & \overline{\chi} \langle -2p \rangle & \overline{\chi}^* \langle 4-2p \rangle & \overline{\chi}        \\
& & & ...  & &;
}$$
the structure of $\HHTct^+= \ffHH(\HTct^+)$ as a $k$-graded $\chi$-$\chi$-bimodule
$$
\xymatrix@C=6pt@R=3pt{
&&&&&& & ... & & \\
&&&&&&           & \overline{\chi}[p-2] & \overline{\chi}^*[p-1]     & \overline{\chi}[3p-4]     & \Omega^{0}[3p-3] \\
&&&&&&           & \overline{\chi}^*[0] & \overline{\chi}[2p-3]          & \Omega^{0}[2p-2]& \\
&&&&&&   & \overline{\chi}[p-2] & \Omega^{0}[p-1] & &\\
&&&&&&  & \Omega^{0}[0]  & & &\\
&&&&&& \underline{\chi}[p-1] & & &&\\
&&&&& \chi &&&   &&  ;   \\
}$$
and finally the structure of $\HHTct^+$ as a $j$-graded $\chi$-$\chi$-bimodule
$$
\xymatrix@C=6pt@R=3pt{
&&&&&& & ... & & \\
&&&&&& & \overline{\chi} \langle p \rangle & \overline{\chi}^* \langle 4+p \rangle & \overline{\chi} \langle 3p\rangle    & \Omega^{0}   \langle 2+3p \rangle \\
&&&&&&           & \overline{\chi}^* \langle 4 \rangle & \overline{\chi}   \langle 2p \rangle       & \Omega^{0}   \langle 2+2p \rangle &\\
&&&&&&   & \overline{\chi} \langle p \rangle & \Omega^{0}  \langle 2+p \rangle && \\
&&&&&&    & \Omega^{0}   \langle 2 \rangle & & &\\
&&&&&& \underline{\chi}   \langle p \rangle & & &&\\
&&&&& \chi \langle 0 \rangle &&&&&    .    \\
}$$

\subsection{Multiplication.}

In order to give the multiplication on $\HHTct$, which thanks to Propositions \ref{cup} and \ref{multstruc} is induced by multiplication in $\Lambda$, 
we first define a number of $\chi$-$\chi$-bimodule homomorphisms between the various components of $\HHTct$.

\begin{lem} \label{crazysymbols}
Let $\bigstar$, $\lozenge_l$, $\lozenge_r$, $\blacklozenge$, $\blacklozenge_r$, 
$\Box_l$, $\Box_r$ and $\blacktriangle$ be the $\chi$-$\chi$-bimodule homomorphisms obtained by applying 
$\HH(\Omega,-)$ to $a:\Theta^\sigma \otimes \Theta^\sigma \rightarrow \Theta$, $\theta_l$, $\theta_r$, $\iota_l$, $\iota_r$, $\nu_l$, $\nu_r$, and $\beta$ from Lemma \ref{homs} respectively, 
which we identify with products of components of $\bbH(\bfc^{\op} \otimes \HTct)$. 
Then the products of basis elements in these spaces that are nonzero are given as follows:

\begin{equation*}
\begin{split}
\bigstar: \quad&\overline{\chi}^* \otimes_\chi \overline{\chi}^* \to \overline{\chi} \\
& \mu_{\frac{p-1}{2}} \otimes \mu_{\frac{p-1}{2}} \mapsto \xi\eta(e_{\frac{p-1}{2}}-e_{\frac{p+1}{2}})\\
&(e_{\frac{p+1}{2}}\otimes e_{\frac{p+1}{2}}x^{\frac{p-3}{2}}y^{\frac{p-1}{2}}e_{\frac{p-1}{2}}) \otimes \mu_{\frac{p-1}{2}}
\mapsto \kappa z^{\frac{p-3}{2}} \\
&(e_{\frac{p-1}{2}}\otimes e_{\frac{p-1}{2}} x^{\frac{p-1}{2}}y^{\frac{p-3}{2}} e_{\frac{p+1}{2}}) \otimes \mu_{\frac{p-1}{2}}
\mapsto  \kappa z^{\frac{p-3}{2}} \\
&\mu_{\frac{p-1}{2}} \otimes (e_{\frac{p+1}{2}}\otimes e_{\frac{p+1}{2}}x^{\frac{p-3}{2}}y^{\frac{p-1}{2}}e_{\frac{p-1}{2}}) 
\mapsto \kappa z^{\frac{p-3}{2}} \\
&\mu_{\frac{p-1}{2}} \otimes(e_{\frac{p-1}{2}}\otimes e_{\frac{p-1}{2}} x^{\frac{p-1}{2}}y^{\frac{p-3}{2}} e_{\frac{p+1}{2}}) 
\mapsto \kappa z^{\frac{p-3}{2}}
\end{split}
\end{equation*}

$$\xymatrix@C=6pt@R=3pt{
\lozenge_l: & \chi \otimes_\chi   \Omega^{0} \to \underline{\chi}, & \quad \lozenge_r: & \Omega^0 \otimes_\chi \chi \to \underline{\chi} \\ 
& 1 \otimes e_p \mapsto z^{p-1} & & e_p \otimes 1 \mapsto z^{p-1}
}$$

$$\xymatrix@C=6pt@R=3pt{
\blacklozenge_l: & \chi \otimes_\chi   \Omega^{0} \to \chi, & \quad \blacklozenge_r: & \Omega^0 \otimes_\chi \chi \to \chi \\ 
& 1 \otimes e_p \mapsto z^{p-1} & & e_p \otimes 1 \mapsto z^{p-1}
}$$

\begin{equation*}
\begin{split}
\Box_l : \quad& \overline{\chi} \otimes_\chi \overline{\chi}^* \to \Omega^{0}, \\
&  1 \otimes (e_s \otimes e_s x^{p-s-1}y^{s-1}) \mapsto e_s, \quad 1 \leq s \leq p-1
\end{split}
\end{equation*}

\begin{equation*}
\begin{split}
\Box_r : \quad& \overline{\chi}^* \otimes_\chi \overline{\chi} \to \Omega^{0}, \\
&  (e_s \otimes e_s x^{p-s-1}y^{s-1}) \otimes 1 \mapsto e_s, \quad 1 \leq s \leq p-1
\end{split}
\end{equation*}

\begin{equation*}
\begin{split} 
\blacktriangle: \quad& \underline{\chi} \otimes_\chi\underline{\chi} \to \Omega^{0},\\ 
& z^{\frac{p-1}{2}} \otimes z^{\frac{p-1}{2}} \mapsto \sum_{s=\frac{p+1}{2}}^{p}e_s
\end{split}
\end{equation*}

\end{lem}
\proof
The product $\bigstar$. Let us consider the element $\kappa z^{\frac{p-3}{2}}$ of $\HH(\Omega, \Theta)$. 
From the proof of Lemma \ref{bimodulecomputation}\eqref{11ii} we find it is equal to 
$\sum_{s=1}^{p-1} a_{s,0} z^{\frac{p-3}{2}}$. 
We know that $a_{s,0}z^{\frac{p-3}{2}}$ is zero unless $s = \frac{p-1}{2}$; consequently
$$\kappa z^{\frac{p-3}{2}} = 
e_{\frac{p-1}{2}} \xi e_{\frac{p+1}{2}} \otimes e_{\frac{p+1}{2}} x^{\frac{p-3}{2}}y^{\frac{p-1}{2}}e_{\frac{p-1}{2}}.$$
The image of $e_{\frac{p-1}{2}} \otimes z^{\frac{p-3}{2}}$ under the differential is 
$$e_{\frac{p-1}{2}} \xi e_{\frac{p+1}{2}} \otimes e_{\frac{p+1}{2}} x^{\frac{p-3}{2}}y^{\frac{p-1}{2}}e_{\frac{p-1}{2}}
- e_{\frac{p+1}{2}} \eta e_{\frac{p-1}{2}} \otimes e_{\frac{p-1}{2}} x^{\frac{p-1}{2}}y^{\frac{p-3}{2}}e_{\frac{p+1}{2}},$$
and therefore in homology we obtain
$$\kappa z^{\frac{p-3}{2}} = 
e_{\frac{p-1}{2}} \xi e_{\frac{p+1}{2}} \otimes e_{\frac{p+1}{2}} x^{\frac{p-3}{2}}y^{\frac{p-1}{2}}e_{\frac{p-1}{2}} =
e_{\frac{p+1}{2}} \eta e_{\frac{p-1}{2}} \otimes e_{\frac{p-1}{2}} x^{\frac{p-1}{2}}y^{\frac{p-3}{2}}e_{\frac{p+1}{2}}.$$
We have $\mu_{\frac{p-1}{2}} = e_{\frac{p-1}{2}} \xi e_{\frac{p+1}{2}} \otimes e_{\frac{p+1}{2}}+ 
e_{\frac{p+1}{2}} \eta e_{\frac{p-1}{2}} \otimes e_{\frac{p-1}{2}}$.
Multiplying in $\bfc^{\op} \otimes \HTct$ gives us $\bigstar$.

The product $\lozenge_l$. Consider the product $\theta_l: \Omega \otimes \Omega^* \rightarrow \Omega e_p \Omega$. 
This factors over the action map $\Omega \otimes \Omega e_p \Omega \rightarrow \Omega e_p \Omega$, and consequently $\lozenge_l$
factors over the action map $\chi \otimes \underline{\chi} \rightarrow \underline{\chi}$.
If we want to know $\lozenge_l$ it therefore suffices to know $\HH(\gamma): \Omega^0 \rightarrow \underline{\chi}$ where again $\gamma$ is the natural surjection $\Omega^* \twoheadrightarrow   \Omega e_p \Omega$ from \eqref{dualseq}.

For every $1\leq s\leq p-1$, the linear form $e_s^*$ vanishes on $\Omega e_p \Omega$ whereas the restriction of $e_p^*$ to $\Omega e_p \Omega$ is equal to $\langle z^{p-1},- \rangle$ (where $\langle -,- \rangle$ is the bilinear form induced by \eqref{omegaiso}). Accordingly, the mapping 
 $\HH(\gamma)$  vanishes on $e_s$ and maps $e_p$ to $z^{p-1}$, which fits with the stated structure of $\lozenge_l$. 

The product $\blacklozenge_l$. The product $\iota_l$ is merely the composition of $\theta_l$ and the embedding of $\Omega e_p \Omega$ in $\Omega$.
Therefore $\blacklozenge_l$ is the composition of $\lozenge_l$ and the natural embedding of $\underline{\chi}$ in $\chi$.

The product $\Box_l$. Consider the product $\nu_l: \Theta \otimes \Theta^\sigma \rightarrow \Omega^*$. 
This is the composite of the action of $\Theta$ on $\Theta^\sigma$ and the embedding $\mu$ of $\Theta^\sigma$ in $\Omega^*$ using \eqref{thetaiso},
in which the socle of $\Theta^\sigma$ is identified with the socle of $\Omega^*$. 
To know $\HH(\nu_l)$ it therefore suffices to know $\HH(\mu)$. 
Since in our computation of $\HH(\Omega, \Omega^*)$ the space $\Omega^0$ is identified with the socle of $\Omega^*$ 
in the tensor product $\bfc^{\op} \otimes \Omega^*$,
and $\mu$ identifies $e_s \otimes e_s x^{p-s-1}y^{s-1}$ with the element of the socle of $\Omega^*$ corresponding to $e_s \in \Omega^0$,
the product $\Box_l$ is as stated.

The products $\lozenge_r$, $\blacklozenge_r$, and $\Box_r$ are established 
similarly to $\lozenge_l$, $\blacklozenge_l$, and $\Box_l$.

The product $\blacktriangle$. We know that under $\blacktriangle$ the radical of $\underline{\chi}$
must have product zero with all elements since $\Omega^0$ is semisimple.
This leaves us with the problem of finding the 
square of the element $z^{\frac{p-1}{2}}$ of $\underline{\chi}$ in $\Omega^0$.
We need to find the element in $\Omega^0$ corresponding to $\beta(z^{\frac{p-1}{2}}\otimes z^{\frac{p-1}{2}})$, that is $\sum_{s=1}^p \beta(z^{\frac{p-1}{2}}\otimes z^{\frac{p-1}{2}})(e_s)e_s$. Now, by the explicit isomorphism described after \eqref{tiltingsquare}  $ \beta(z^{\frac{p-1}{2}}\otimes z^{\frac{p-1}{2}})(e_s) = \langle e_s z^{\frac{p-1}{2}}e_s, e_s z^{\frac{p-1}{2}}e_s\rangle$, which equals $1$ if $s\geq \frac{p+1}{2}$ and $0$ otherwise. Thus the resulting element in $\Omega^0$ is $\sum_{s=\frac{p+1}{2}}^{p}e_s$, as stated.
\endproof

We use these maps to describe the product in $\HHTct$, where we again gather together components which are isomorphic (up to shift), according to whether they lie in $\HHTct^+$ or $\HHTct^-$, in a similar way as in Proposition \ref{clubsuitmult}.

\begin{thm} \label{spadesuitproduct}
Products  between the various components in $\HHTct$ are given by the following table
$$
\xymatrix@C=10pt@R=3pt{ & \chi_- & \overline{\chi}_- & \overline{\chi}_-^* & \underline{\chi} & \overline{\chi}_+ & \overline{\chi}_+^* & \Omega^0_+\\
\chi_- &   a &a&a&a&0&0&     \blacklozenge, \lozenge, a        \\
\overline{\chi}_- & a&a&a&0& 0,a & 0, a, \Box &0\\
\overline{\chi}_-^* & a&a&\bigstar&0&0,a,\Box &0, \bigstar &0  \\
\underline{\chi} & a&0&0& \blacktriangle & 0&0&0\\
\overline{\chi}_+ & 0 &0,a&0,a,\Box &0&0&0&0\\
\overline{\chi}_+^* & 0 & 0,a,\Box & 0, \bigstar &0&0&0&0 \\
 \Omega^0_+ &  \blacklozenge, \lozenge, a & 0&0&0&0&0&0
}$$
Possible ambiguities are covered by further tables. For the product of $\Omega^0_+$ and $\chi_-$:
$$\xymatrix@C=10pt@R=3pt{
\textrm{Component in which the product lands:} & \chi & \underline{\chi} & \Omega^0_+ \\
\textrm{Natural map describing the product:} & \blacklozenge & \lozenge & a
}$$
For the product of $\overline{\chi}_+$ and $\overline{\chi}_-$:
$$\xymatrix@C=10pt@R=3pt{
\textrm{Component in which the product lands:} & \overline{\chi}_+ & \HHTct^- & \\
\textrm{Natural map describing the product:} & a & 0
}$$
For the product of ${\overline{\chi}^*}_-$ and $\overline{\chi}_+$:
$$\xymatrix@C=10pt@R=3pt{
\textrm{Component in which the product lands:} &{\overline{\chi}^*}_+ & \underline{\chi}_+ & \Omega^0_+ & \HHTct^- \\
\textrm{Natural map describing the product:} & a & 0 & \Box & 0
}$$
For the product of ${\overline{\chi}^*}_+$ and $\overline{\chi}_-$:
$$\xymatrix@C=10pt@R=3pt{
\textrm{Component in which the product lands:} & {\overline{\chi}^*}_+ & \underline{\chi}_+ & \Omega^0_+ & \HHTct^- \\
\textrm{Natural map describing the product:} & a & 0 & \Box & 0
}$$
For the product of $\overline{\chi}^*_+$ and $\overline{\chi}^*_-$:
$$\xymatrix@C=10pt@R=3pt{
\textrm{Component in which the product lands:} & \overline{\chi}_+ & \HHTct^- & \\
\textrm{Natural map describing the product:} & \bigstar & 0
}$$

\end{thm}
\proof
All the action products are inherited from action products in $\HTct$; all other nonzero products are inherited from nonzero products in $\HTct$ or via Lemma \ref{crazysymbols}. The zero products  are either inherited from zero products in $\HTct$,
or determined by the fact that the products lie in degrees in which there are no nonzero elements with respect to the various gradings;
for example $\HH(\epsilon) = \HH(\zeta) = 0$ by this reasoning.
\endproof

\section{A monomial basis.} \label{final}

As any Ringel self-dual block of polynomial representations of $G$ is equivalent to $\bbO_{F,0} \bbO_{\bfc,\ut}^l(F,(F,F))\ml$ for some $l \geq 0$, we have established the following:

\begin{thm}
We have isomorphisms of $k$-graded algebras
$$\hh_l \cong \fO_F \fO_{\HHTct}^l(F[z,z^{-1}]).$$
\end{thm}
\proof This is a restatement of Proposition \ref{stock}.
\endproof

We describe a basis for $\HHTct$ indexed by elements of a polytope.
Roughly, we label basis elements $m_{d,e}$ for $\HHTct$ by a pair $(d,e)$ where $d \in \mathbb{Z}^3$
denotes a $ijk$-degree, 
and $e$ denotes an element of $\Omega^0$, either $1$ or an idempotent.

More precisely, here is our basis for $\HHTct$:
\begin{equation*}\begin{split}
\mathbf{B}_{\HHTct} = & \mathbf{B}_{\chi_-} \cup \mathbf{B}_{\overline{\chi}_-} \cup \mathbf{B}_{{\overline{\chi}^*}_-} \cup \mathbf{B}_{\underline{\chi}}  \cup
\mathbf{B}_{\overline{\chi}_+}  \cup \mathbf{B}_{{\overline{\chi}^*}_+} \cup \mathbf{B}_{\Omega^0} \\
= & \{ m_{a,b,i,j+ap,k+a(1-p),e} | m_{j,k,e} \in \mathcal{B}_\chi, a \leq 0, b = 0, i = a+b \}\\
\cup & \{ m_{a,b,i,j+ ap,k+a(1-p),e} | m_{j,k,e} \in \mathcal{B}_{\overline{\chi}}, a \leq 0, b \leq -2, b \textrm{ even }, i = a+b  \}\\
\cup & \{ m_{a,b,i,j+(4-p)+ap,k+(p-2)+a(1-p),e} | m_{j,k,e} \in \mathcal{B}_{{\overline{\chi}^*}}, a \leq 0, b \leq -1, b \textrm{ odd }, i = a+b \}\\
\cup & \{ m_{1,0,1,j+p,k+1-p,e} | m_{j,k,e} \in \mathcal{B}_{\underline{\chi}}  \} \\
\cup & \{ m_{a,b,i,j+(a-1)p,k+1+(a-1)(1-p),e} | m_{j,k,e} \in \mathcal{B}_{\overline{\chi}}, a \geq 2, b \geq 1, b \textrm{ odd }, i = a+b  \}\\
\cup & \{ m_{a,b,i,j+4+(a-2)p,k+(a-2)(1-p),e} | m_{j,k,e} \in \mathcal{B}_{{\overline{\chi}^*}}, a \geq 2, b \geq 2, b \textrm{ even }, i = a+b  \}\\
\cup & \{ m_{a,b,i,j+2+(a-2)p,k+(a-2)(1-p),e} | m_{j,k,e} \in \mathcal{B}_{\Omega^0}, a \geq 2, b=0, i = a+b  \}
\end{split}\end{equation*}
We describe the $a,b$ grading as follows: 
in our pictures of $\HHTct$ a shift by $a$ corresponds to a move to the northeast by $a$ and a shift by $b$ corresponds 
to a move to the north by $b$.
The product of a pair of basis elements in $\HHTct$ is either another basis element, 
or the sum of a basis element and the negative of another basis element, or $\pm \frac{1}{2}$ a basis element, or zero; 
when a product of $m_{a,b,i,j,k,e} . m_{a',b',i',j',k',e'}$ is nonzero, 
the basis elements in the product take the form $m_{a+a',b+b',i+i',j+j',k+k',y}$.  
Precise formulas for the product are given by the formulas in the statement of Lemma \ref{crazysymbols} 
and the table in the statement of Theorem \ref{spadesuitproduct}.

We can now use this to construct a basis for $\hh_l$.

\begin{cor} \label{totalbasis}
The algebra $\hh_l$ inherits an explicit basis from $\HHTct$.
\end{cor}

Before proving this, we recall that the $ik$-homogeneous component of $\fO_{\HHTct}^l (F[z,z^{-1}])$ is given by
$\bigoplus\HHTct^{ij_1k_1}\otimes \HHTct^{j_1j_2k_2}\otimes \cdots \HHTct^{j_{l-1}j_1k_l}\otimes z^{k_l}$
where the sum runs over all integers $j_1,\dots, j_l$ and $k_1,\dots, k_l$ such that $k_1+\cdots +k_l=k$. The operator $\fO_F$ then projects onto the homogeneous component of $i$-degree $0$.

\proof
We explicitly write down such a basis as follows: let $\mathbf{B}_{\HHTct}$ denote our basis for $\HHTct$.
We have a basis for the algebra $\HHTct^{\otimes_F l} \otimes_F F[z,z^{-1}]$ given by
$\mathbf{B}_{\HHTct}^{\times l} \times \{z^d | d \in \mathbb{Z}\}$;
the product of basis elements is the super $\times$ product.
We define the \emph{weight} of a monomial
$m_{w^1} \otimes ... \otimes m_{w^q} \otimes z^\alpha$ in $\mathbf{B}_{\HHTct}^{\times l} \times \{z^d | d \in \mathbb{Z}\}$ to be
$$(w^2_i - w^1_j,w^3_i - w^2_j,....,w^l_i - w^{l-1}_j, \alpha - w^l_j) \in \mathbb{Z}^{l+1},$$
where $(w_i,w_j)$ denotes the $ij$-degree of $m_w$.
We then have a basis for the algebra $\fO_F \fO_{\HHTct}^l (F[z,z^{-1}])$ given by weight zero elements in
$\mathbf{B}_{\HHTct}^{\times l} \times \{z^d | d \in \mathbb{Z}\}$;
the product is the restriction of the product on $\mathbf{B}_{\HHTct}^{\times l} \times \{z^d | d \in \mathbb{Z}\}$.
\endproof

\begin{cor}
The map $\hh_l \rightarrow \hh_{l-1}$ is surjective for $l \geq 1$.
\end{cor}
\proof The map $\HHTct \rightarrow F$ is surjective,
implying $$\fO_{\HHTct}(a) \rightarrow \fO_F(a)$$ is surjective for any $a$,
implying $$\fO_F \fO_{\HHTct}(a) \rightarrow \fO_F^2(a) = \fO_F(a)$$ is surjective for any $a$,
implying $$\fO_F \fO_{\HHTct}^l(F[z,z^{-1}]) \rightarrow \fO_F\fO_{\HHTct}^{l-1}(F[z,z^{-1}])$$ is surjective,
implying $\hh_l \rightarrow \hh_{l-1}$ is surjective.
\endproof

\normalfont

{\sc Vanessa Miemietz}\\
School of Mathematics, University of East Anglia, Norwich, NR4 7TJ, UK, \\{\tt v.miemietz@uea.ac.uk}\\
{\sc Will Turner}\\
Department of Mathematics, University of Aberdeen, Fraser Noble Building, King's College, Aberdeen AB24 3UE, UK, {\tt w.turner@abdn.ac.uk}.

\begin{thebibliography}{00}
\normalfont

\bibitem{BGS} A. Beilinson, V. Ginzburg, W. Soergel, \emph{Koszul duality patterns in representation theory}, J. Amer. Math. Soc.  9  (1996),  no. 2, 473--527.

\bibitem{BGMS} R. Buchweitz, E. Green, D. Madsen, and \O. Solberg, \emph{Finite Hochschild cohomology without finite global
dimension}, Math. Res. Lett. 12 (2005), 805--816.

\bibitem{CPS}  E. Cline, B. Parshall and L. Scott, \emph{Finite-dimensional algebras and highest weight categories.}
J. Reine Angew. Math. 391 (1988), 85--99.


\bibitem{DR} V. Dlab, C.M. Ringel, \emph{A construction for quasi-hereditary algebras}, Compositio Math., 70 (1989) no.~2, 155--175.
\bibitem{Ger} M. Gerstenhaber, \emph{The cohomology structure of an associative ring}.  Ann. of Math. (2) 78 (1963) 267--288.

\bibitem{Green} J. A. Green, \emph{Polynomial representations of ${\rm
GL}_{n}$}, Lecture Notes in Mathematics, 830. Springer, Berlin,
1980.

\bibitem{EH} K. Erdmann and A. Henke, \emph{On Ringel duality for Schur algebras}, Math. Proc.
Cambridge Philos. Soc. 132(1) (2002), 97--116.

\bibitem{Ha} D. Happel, \emph{Triangulated categories in the representation theory of finite-dimensional algebras}, London Mathematical Society Lecture Notes 119, Cambridge University Press,1988.

\bibitem{Ke1}
B. Keller, \emph{Deriving dg categories},  Ann. Sci. \'Ecole Norm. Sup. (4) 27 (1994), no.~1, 63--102.

\bibitem{Ke2} B. Keller, \emph{$A_\infty$-algebras, modules and functor categories}, Trends in representation theory of algebras and related topics, Contemp. Math. 406, 67--93. Amer. Math. Soc., Providence, RI, 2006.

\bibitem{Ke}
B. Keller, \emph{On differential graded categories}, International Congress of
  Mathematicians. Vol. II, Eur. Math. Soc., Z\"urich, 2006, pp.~151--190.



\bibitem{Kr} U. Kr\"ahmer, \emph{Notes on Koszul algebras}, http://www.maths.gla.ac.uk/~ukraehmer/connected.pdf.

\bibitem{LH} K. Lef\`evre-Hasegawa, Sur les $A_\infty$ cat\'egories, PhD thesis, Universit\'e Paris 7 - Denis Diderot, 2003.

\bibitem{Maz} V. Mazorchuk, \emph{Koszul duality for stratified algebras, I}, Balanced quasi-hereditary algebras. Manuscripta Math. 131 (2010), no. 1--2, 1--10.

\bibitem{MT2} V. Miemietz, W. Turner, \emph{Homotopy, Homology and $GL_2$}, Proc. London Math. Soc. (3) 100 (2010), no.~2, 585--606.

\bibitem{MT3} V. Miemietz, W. Turner, \emph{Koszul dual $2$-functors and extension algebras of simple modules for $GL_2$}, 
Selecta Math. (N.S.) 21 (2015), no.~2,  605--648.

\bibitem{MT4} V. Miemietz, W. Turner, \emph{The Weyl extension algebra of $GL_2(\bar{\mathbb{F}}_p)$},  Adv. Math.  246 (2013), 144--197.

\bibitem{Ne} C. Negron, \emph{The cup product on Hochschil cohomology via twisting cochains and applications to Koszul rings}, J. of Pure and Applied Algebra  221 (2017), 1112--1133

\bibitem{Rickard}  J. Rickard, \emph{Derived equivalences as derived functors}, J. London Math. Soc. (2) 43 (1991), no. 1, 37--48.

\bibitem{RR} R. Rouquier, \emph{Derived equivalences and finite
dimensional algebras}, Proceedings of the International Congress of
Mathematicians (Madrid, 2006 ), vol II, pp. 191-221, EMS Publishing
House, 2006.

\bibitem{Salfelder} F. Salfelder, \emph{Hochschild cohomology of category $\mathcal{O}$}, felix.salfelder.org/misc/HH.ps
\bibitem{Snashall} N. Snashall,\emph{Support varieties and the Hochschild cohomology ring modulo nilpotence}, Proceedings of the 41st Symposium on Ring Theory and Representation Theory, 68--82, Symp. Ring Theory Represent. Theory Organ. Comm., Tsukuba, 2009.


\bibitem{Xu} F. Xu, \emph{Hochschild and ordinary cohomology rings of small categories}, Adv. Math. 219 (2008), 1872--1893.





\end{thebibliography}
\end{document}